\documentclass[11pt,reqno,oneside]{article}

\usepackage{multirow}
\usepackage{lineno}
\usepackage{pstricks,pst-node,pst-tree}
\usepackage{amsmath}
\usepackage{amssymb}
\usepackage[lmargin=.55in, rmargin=.55in, tmargin=1in, bmargin=1in]{geometry}
\usepackage{xcolor}
\usepackage{graphicx}
\usepackage[font=small]{caption}
\usepackage{natbib}
\usepackage[labelformat = empty,position=top]{subcaption}
\usepackage{pstricks}

\newtheorem{lemm}{Lemma}
\newtheorem{prop}{Proposition}

\definecolor{myblue}{rgb}{0,0.1,0.9}
\definecolor{myred}{rgb}{0.9,0.1,0}

\title{{\bf The distributions under two species-tree models \\ of the number of root ancestral configurations \\ for matching gene trees and species trees}}

\author{Filippo Disanto\thanks{Department of Mathematics, University of Pisa, Pisa 56126, Italy. Email: filippo.disanto@unipi.it.}, Michael Fuchs\thanks{Department of Mathematical Sciences, National Chengchi University, Taipei 116, Taiwan. Email: mfuchs@nctu.edu.tw.}, Ariel R.~Paningbatan\thanks{Institute of Mathematics, University of the Philippines Diliman, Quezon City 1101, Philippines. Email: arpaningbatan@math.upd.edu.ph.}, Noah A.~Rosenberg\thanks{Department of Biology, Stanford University, Stanford, CA 94305, USA. Email: noahr@stanford.edu.}}

\begin{document}

\maketitle

\begin{abstract}
For a pair consisting of a gene tree and a species tree, the \emph{ancestral configurations} at an internal node of the species tree are the distinct sets of gene lineages that can be present at that node. Ancestral configurations appear in computations of gene tree probabilities under evolutionary models conditional on fixed species trees, and the enumeration of \emph{root} ancestral configurations---ancestral configurations at the root of the species tree---assists in describing the complexity of these computations. In the case that the gene tree matches the species tree in topology, we study the distribution of the number of root ancestral configurations of a random labeled tree topology under each of two models. First, choosing a tree uniformly at random from the set of labeled topologies with $n$ leaves, we extend an earlier computation of the asymptotic exponential growth of the mean and variance of the number of root ancestral configurations, showing that the number of root ancestral configurations of a random tree asymptotically follows a lognormal distribution; the logarithm has mean $\sim$0.272$n$ and variance $\sim$0.034$n$. The asymptotic mean of the logarithm of the number of root ancestral configurations produces $e^{0.272n} \approx 1.313^n$ when exponentiated, numerically close to the previously obtained mean of $(4/3)^n$ for the exponential growth of the number of root ancestral configurations. Next, considering labeled topologies selected according to the Yule--Harding model, we obtain the asymptotic mean and variance of the number of root ancestral configurations of a random tree and the asymptotic distribution of its logarithm. The asymptotic mean follows $\sim$1.425$^n$ and the variance follows $\sim$2.045$^n$; the random variable has an asymptotic lognormal distribution, and its logarithm has mean $\sim$0.351$n$ and variance $\sim$0.008$n$. The asymptotic mean of the logarithm produces $e^{0.351n} \approx 1.420^n$ when exponentiated, close to the mean of $1.425^n$. With the higher probabilities assigned by the Yule--Harding model to balanced trees in comparison with those assigned under the uniform model, a larger asymptotic exponential growth $\sim$1.425$^n$ of the mean number of root ancestral configurations for the Yule--Harding model compared to $(4/3)^n$ in the uniform model suggests an effect of increasing tree balance in increasing the number of root ancestral configurations.
\end{abstract}

\begin{itemize}
\item {\bf Keywords:} analytic combinatorics, gene trees, phylogenetics, species trees.
\item {\bf Mathematics subject classification (2010):} 05A15 $\cdot$ 05A16 $\cdot$ 92B10 $\cdot$ 92D15
\item {\bf Running title:} Root ancestral configurations

% \item {\bf Possible target journals:}
% math journals such as \emph{Ann Combinator}, \emph{Discr Appl Math}, \emph{SIAM J Discr Math};
% mathematical biology journals such as \emph{Bull Math Biol}, \emph{IEEE/ACM Trans Comp Biol Bioinf}, \emph{J Comput Biol}, \emph{J Math Biol}.
\end{itemize}

\clearpage

%%%%%%%%%%%%%%%%%%%%%%%%%%%%%%%%%%%%%%%%%%%%%%%%%%%%%%%%%%%%%%
%%%%%%%%%%%%%%%%%%%%%%%%%%%%%%%%%%%%%%%%%%%%%%%%%%%%%%%%%%%%%%
%%%%%%%%%%%%%%%%%%%%%%%%%%%%%%%%%%%%%%%%%%%%%%%%%%%%%%%%%%%%%%

\section{Introduction}

In the study of combinatorial properties of species trees, trees that describe evolutionary relationships among species, and gene trees, trees that describe evolutionary relationships among gene lineages for members of the species, one useful concept is that of an ancestral configuration \citep{Wu12, DisantoAndRosenberg17}. Given a gene tree, a species tree, and a node of the species tree, an ancestral configuration is a list of the gene lineages that are present at the node of the species tree (Figure 1). Looking backward in time, or from the tips of trees to the root, the fact that gene lineages only find their common ancestors once their associated species have found common ancestors produces conditions describing which ancestral configurations are present at a species tree node. These conditions enable the enumeration of the configurations. Ancestral configurations appear in recursive evaluations of the probabilities of gene tree topologies conditional on species tree topologies \citep{Wu12}, so that enumerations of ancestral configurations assist in assessing the complexity of the computation.

When the node at which an ancestral configuration is considered is the root node of the species tree, ancestral configurations are termed \emph{root ancestral configurations}, or root configurations for short. For matching gene trees and species trees---that is, if the species tree and gene tree have the same labeled topology---the number of root configurations is greater than or equal to the number of ancestral configurations for any other species tree node. This property can be used to show that as the number of taxa increases, the total number of ancestral configurations for the gene tree and species tree---the sum of the number of ancestral configurations across all species tree nodes---has the same exponential growth as the number of root configurations \citep[][Section 2.3.2]{DisantoAndRosenberg17}. Hence, it suffices for investigations of the exponential growth of the total number of ancestral configurations for matching gene trees and species trees to focus on root configurations.

\cite{DisantoAndRosenberg17} studied the number of root configurations for matching gene trees and species trees, considering the number of root configurations of families of increasingly large trees. They characterized the labeled tree topologies with the largest number of root configurations among trees with $n$ leaves, showing that this number of root configurations lies between $k_0^{n-1/4}-1$ and $k_0^n-1$, where $k_0$ is a constant approximately equal to 1.5028 \citep[][Proposition 4]{DisantoAndRosenberg17}. They then studied the number of root configurations in trees selected uniformly at random from the set of labeled topologies with $n$ leaves. Using techniques of analytic combinatorics, they showed that the mean number of root configurations grows with $(4/3)^n$, and the variance with $\sim$1.8215$^n$ \citep[][Propositions 5 and 6]{DisantoAndRosenberg17}.

Here, we extend these results on the distribution of the number of root configurations under a model imposing a uniform distribution on the set of labeled topologies. We obtain an asymptotic normal distribution for the logarithm of the number of root configurations under the uniform model, finding that its mean, approximately $0.272n$, generates exponential growth $e^{0.272n} \approx 1.313^n$. We next obtain similar results under the Yule--Harding model, including the asymptotic mean and variance of the number of root configurations and the asymptotic distribution of its logarithm.

%%%%%%%%%%%%%%%%%%%%%%%%%%%%%%%%%%%%%%%%%%%%%%%%%%%%%%%%%%%%%%
%%%%%%%%%%%%%%%%%%%%%%%%%%%%%%%%%%%%%%%%%%%%%%%%%%%%%%%%%%%%%%
%%%%%%%%%%%%%%%%%%%%%%%%%%%%%%%%%%%%%%%%%%%%%%%%%%%%%%%%%%%%%%

\section{Preliminaries}

We study ancestral configurations for rooted binary leaf-labeled trees. In Section \ref{secTreeClasses}, we introduce results on various classes of trees. In Section~\ref{yhdistr}, we discuss the Yule--Harding distribution on labeled topologies. In Section~\ref{analyticomb}, we recall properties of generating functions and analytic combinatorics. Following \cite{Wu12}, in Section~\ref{igno} we define ancestral configurations, and we review enumerative results from \cite{DisantoAndRosenberg17}. In Section \ref{addp}, we relate ancestral configurations to the additive tree parameters of \cite{Wagner15}.

%%%%%%%%%%%%%%%%%%%%%%%%%%%%%%%%%%%%%%%%%%%%%%%%%%%%%%%%%%%%%%
%%%%%%%%%%%%%%%%%%%%%%%%%%%%%%%%%%%%%%%%%%%%%%%%%%%%%%%%%%%%%%
%%%%%%%%%%%%%%%%%%%%%%%%%%%%%%%%%%%%%%%%%%%%%%%%%%%%%%%%%%%%%%

\subsection{Classes of trees}
\label{secTreeClasses}

We will need to consider many classes of trees: labeled topologies, unlabeled topologies, ordered unlabeled topologies, labeled histories, unlabeled histories, and ordered unlabeled histories.

\subsubsection{Labeled topologies}
\label{filo}

\begin{figure}
\begin{center}
\includegraphics*[scale=0.76,trim=0 0 0 0]{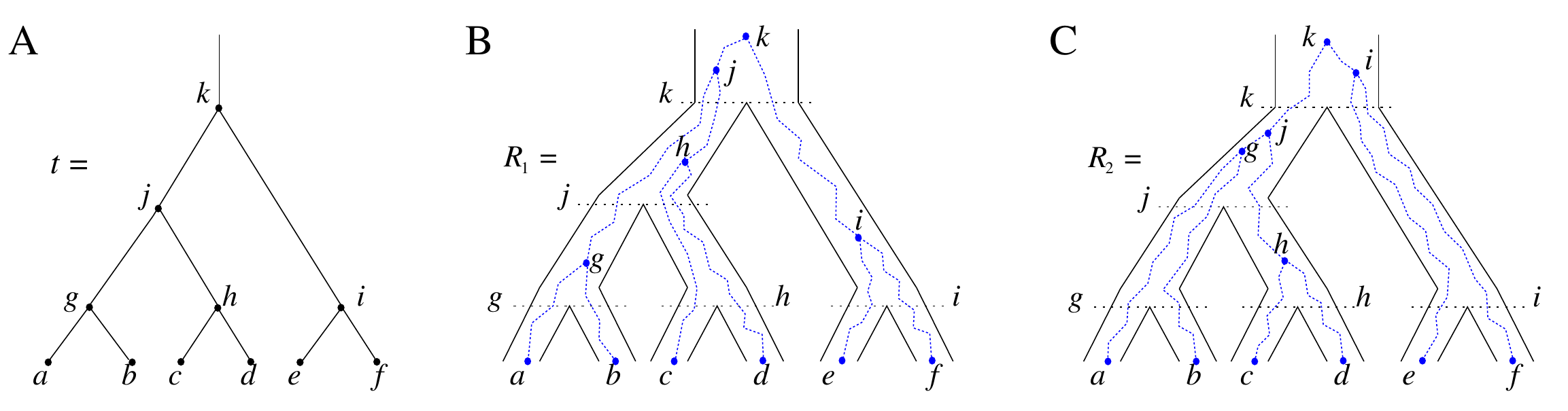}
\end{center}
\vspace{-.7cm}
\caption{{\small A gene tree and species tree with matching labeled topology $t$. {\bf(A)} A tree $t$ of size $6$, characterized by its shape and taxon labels. For convenience, we label the internal nodes of $t$, by $g,h,i,j,k$ in this case, identifying each lineage (edge) by its immediate descendant node. For example, lineage $h$ results from coalescence of lineages $c$ and $d$. {\bf (B)} A possible realization $R_1$ of the gene tree in (A) (dotted lines) in the matching species tree (solid lines). The ancestral configurations at species tree nodes $j$ and $k$ are $\{ g,c,d \}$ and $\{ g,h,i \}$, respectively. {\bf (C)} A different realization $R_2$ of the gene tree in (A) in the species tree. At species tree nodes $j$ and $k$, the configurations are $\{ a,b,h \}$ and $\{ j,e,f \}$, respectively.}}
\label{configa}
\end{figure}

We refer to a bifurcating rooted tree $t$ with $|t|=n$ labeled leaves as a \emph{labeled topology} of size $|t|=n$, or a ``tree'' for short (Fig.~\ref{configa}A); these trees are sometimes called phylogenetic trees or Schr\"{o}der trees. For the set $\{a,b,c, \ldots \}$ of possible labels for the taxa of a tree, we impose an alphabetical linear order $a \prec b \prec c \prec \ldots $ The leaf labels of a tree of size $n$ are the first $n$ labels in the order $\prec$.

We denote by $T_n$ the set of trees of size $n$, with $T = \bigcup_{n=1}^\infty T_n$ denoting the set of all trees. The number of trees of size $n \geq 2$ is $|T_n|=(2n-3)!! = 1 \times 3 \times 5 \times \ldots \times (2n-3)$
\citep{Felsenstein78}, or, for $n \geq 1$,
\begin{equation}
\label{carciofo}
|T_n| =
\frac{(2n-2)!}{2^{n-1} (n-1)!} =
\frac{(2n)!}{2^n(2n-1) n!}.
\end{equation}
The exponential generating function for $|T_n|$ is
\begin{equation*}
%\label{genio}
T(z)= \sum_{t \in T} \frac{z^{|t|}}{|t|!} = \sum_{n=1}^\infty \frac{|T_n| z^n}{n!} = z + \frac{z^2}{2} + \frac{3z^3}{6} + \frac{15z^4}{24} + \ldots,
\end{equation*}
given by \cite[][Example II.19]{FlajoletAndSedgewick09}
\begin{equation}
\label{expo}
T(z) = 1-\sqrt{1-2z}.
\end{equation}
%%%%%%%%%%%%%%%%%%%%%%%%%%%%%%%%%%%%%%%%%%%%%%%%%%%%%%%%%%%%%%
%%%%%%%%%%%%%%%%%%%%%%%% Figure 1 %%%%%%%%%%%%%%%%%%%%%%%%%%%%
%%%%%%%%%%%%%%%%%%%%%%%%%%%%%%%%%%%%%%%%%%%%%%%%%%%%%%%%%%%%%%
%%%%%%%%%%%%%%%%%%%%%%%%%%%%%%%%%%%%%%%%%%%%%%%%%%%%%%%%%%%%%%

%%%%%%%%%%%%%%%%%%%%%%%%%%%%%%%%%%%%%%%%%%%%%%%%%%%%%%%%%%%%%%
%%%%%%%%%%%%%%%%%%%%%%%%%%%%%%%%%%%%%%%%%%%%%%%%%%%%%%%%%%%%%%
%%%%%%%%%%%%%%%%%%%%%%%%%%%%%%%%%%%%%%%%%%%%%%%%%%%%%%%%%%%%%%

\subsubsection{Ordered unlabeled topologies}
\label{out}

\begin{figure}
\begin{center}
\includegraphics*[scale=0.44,trim=0 0 0 0]{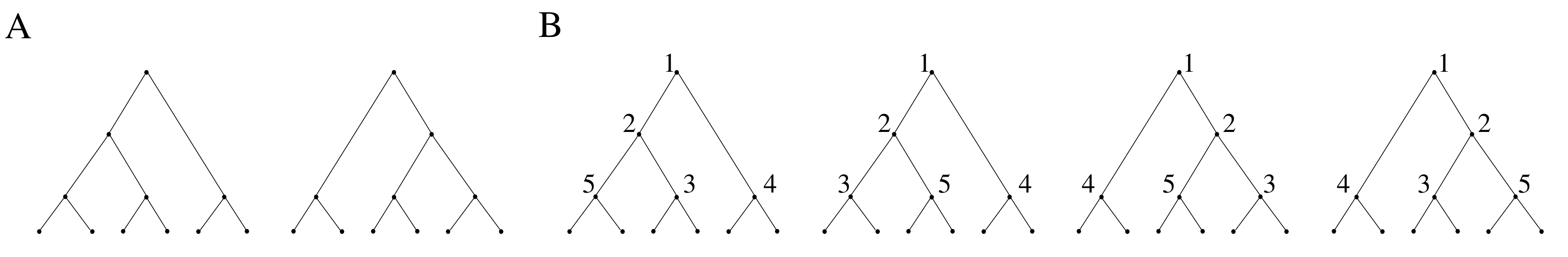}
\end{center}
\vspace{-.7cm}
\caption{{\small
Ordered unlabeled topologies and histories. {\bf(A)} The two orientations of the unlabeled topology that underlies the labeled topology of Fig.~\ref{configa}A. {\bf(B)} The four orientations of the unlabeled history underlying the labeled history in Fig.~\ref{labhist}A.}}
\label{orient}
\end{figure}

An \emph{orientation} of an unlabeled topology $t$ is a planar embedding of $t$ in which subtrees descending from the internal nodes of $t$ are considered with a left--right orientation. For instance, the \emph{unlabeled topology} underlying the labeled topology depicted in Fig.~\ref{configa}A has exactly two different orientations, which are depicted in Fig.~\ref{orient}A. An orientation of an unlabeled topology is called an \emph{ordered} unlabeled topology. The set of all possible ordered unlabeled topologies of size $n$ is enumerated by the Catalan number $C_{n-1}$ \cite[][Exercise 6.19d]{Stanley99}, where
\begin{equation}
\label{cata}
C_n = \frac{1}{n+1} {{2n}\choose{n}}.
\end{equation}
The ordinary generating function is
\begin{equation*}
%\label{catag}
C(z) = \sum_{n=0}^\infty C_n z^n = \frac{1-\sqrt{1-4z}}{2z}.
\end{equation*}
Ordered unlabeled topologies are also called ``pruned binary trees,'' for example by \cite{Wagner15} (see also \cite{FlajoletAndSedgewick09}, Example I.13).

%%%%%%%%%%%%%%%%%%%%%%%%%%%%%%%%%%%%%%%%%%%%%%%%%%%%%%%%%%%%%%
%%%%%%%%%%%%%%%%%%%%%%%%%%%%%%%%%%%%%%%%%%%%%%%%%%%%%%%%%%%%%%
%%%%%%%%%%%%%%%%%%%%%%%%%%%%%%%%%%%%%%%%%%%%%%%%%%%%%%%%%%%%%%

\subsubsection{Labeled histories}

\begin{figure}
\begin{center}
\includegraphics*[scale=0.64,trim=0 0 0 0]{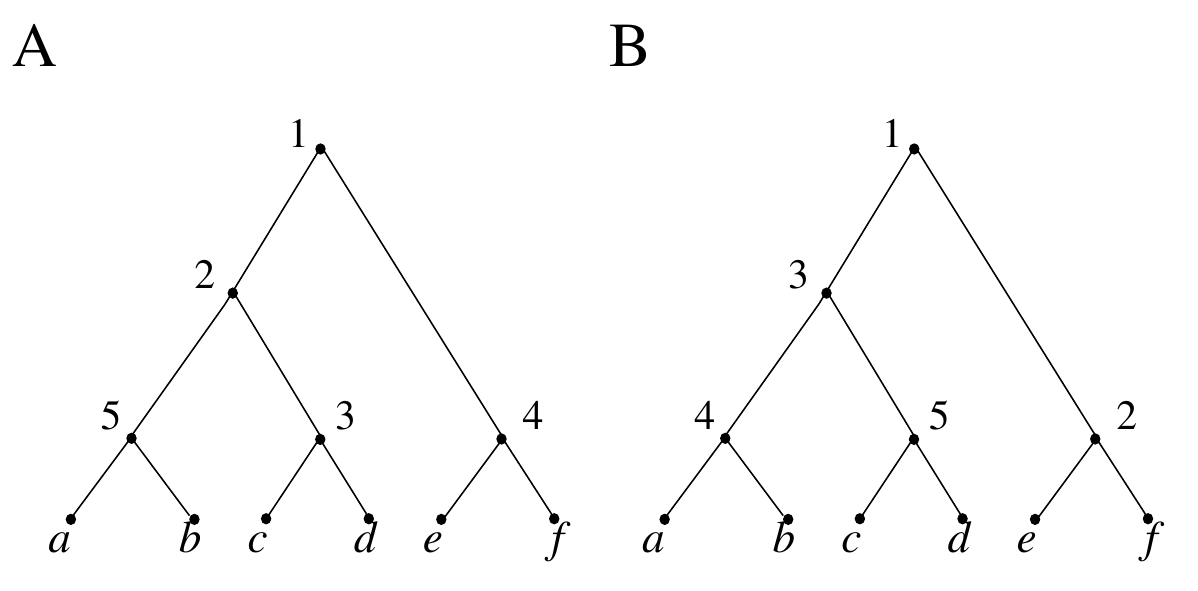}
\end{center}
\vspace{-.7cm}
\caption{{\small
Labeled histories. {\bf(A)} The labeled history of the labeled gene tree topology depicted in Fig.~\ref{configa}B. The temporal ordering of the coalescence events in the gene tree is determined by the integer labeling of the internal nodes of the associated labeled topology. {\bf(B)} The labeled history of the labeled gene tree topology depicted in Fig.~\ref{configa}C. }}
\label{labhist}
\end{figure}

A \emph{labeled history} is a labeled topology together with a temporal (linear) ordering of its internal nodes (Fig.~\ref{labhist}). If $t$ is a labeled history of size $n$, then we represent the time ordering of its $n-1$ bifurcations by bijectively associating each internal node of $t$ with an integer label in the interval $[1, n-1]$. The labeling is increasing in the sense that each internal node other than the root has a larger label than its parent node.

%%%%%%%%%%%%%%%%%%%%%%%%%%%%%%%%%%%%%%%%%%%%%%%%%%%%%%%%%%%%%%
%%%%%%%%%%%%%%%%%%%%%%%% Figure 2 %%%%%%%%%%%%%%%%%%%%%%%%%%%%
%%%%%%%%%%%%%%%%%%%%%%%%%%%%%%%%%%%%%%%%%%%%%%%%%%%%%%%%%%%%%%
%%%%%%%%%%%%%%%%%%%%%%%%%%%%%%%%%%%%%%%%%%%%%%%%%%%%%%%%%%%%%%

%%%%%%%%%%%%%%%%%%%%%%%%%%%%%%%%%%%%%%%%%%%%%%%%%%%%%%%%%%%%%%
%%%%%%%%%%%%%%%%%%%%%%%% Figure 3 %%%%%%%%%%%%%%%%%%%%%%%%%%%%
%%%%%%%%%%%%%%%%%%%%%%%%%%%%%%%%%%%%%%%%%%%%%%%%%%%%%%%%%%%%%%
%%%%%%%%%%%%%%%%%%%%%%%%%%%%%%%%%%%%%%%%%%%%%%%%%%%%%%%%%%%%%%

For a given label set of size $n$, the set of labeled histories is denoted $H_n$. Its cardinality is~\citep[][p.~46]{Steel16}
\begin{equation}
\label{carciofino}
|H_n| = \frac{n! \, (n-1)!}{2^{n-1}}.
\end{equation}

%%%%%%%%%%%%%%%%%%%%%%%%%%%%%%%%%%%%%%%%%%%%%%%%%%%%%%%%%%%%%%
%%%%%%%%%%%%%%%%%%%%%%%%%%%%%%%%%%%%%%%%%%%%%%%%%%%%%%%%%%%%%%
%%%%%%%%%%%%%%%%%%%%%%%%%%%%%%%%%%%%%%%%%%%%%%%%%%%%%%%%%%%%%%

\subsubsection{Ordered unlabeled histories}
\label{ouh}

By removing taxon labels of a labeled history $t$, we obtain the unlabeled history underlying $t$. As we did for unlabeled topologies, we define an orientation of an unlabeled history $t$ as a planar embedding of $t$ in which child nodes are considered with a left--right orientation. Fig.~\ref{orient}B shows the orientations of the unlabeled history underlying the labeled history of Fig.~\ref{labhist}A. We call an orientation of an unlabeled history an \emph{ordered unlabeled history}. The set of all ordered unlabeled histories of size $n$ is enumerated by $F_{n-1}$ \citep[][p.~47]{Steel16}, where
\begin{equation}
\label{fact}
F_n = n!.
\end{equation}
Ordered unlabeled histories are also called ``binary increasing trees'' \citep[]{BergeronEtAl92,Wagner15} or ``ranked oriented trees'' \citep{Steel16}.

%%%%%%%%%%%%%%%%%%%%%%%%%%%%%%%%%%%%%%%%%%%%%%%%%%%%%%%%%%%%%%
%%%%%%%%%%%%%%%%%%%%%%%%%%%%%%%%%%%%%%%%%%%%%%%%%%%%%%%%%%%%%%
%%%%%%%%%%%%%%%%%%%%%%%%%%%%%%%%%%%%%%%%%%%%%%%%%%%%%%%%%%%%%%

\subsection{The Yule--Harding distribution}
\label{yhdistr}

Different labeled histories can share the same underlying labeled topology. For example, the labeled histories of Fig.~\ref{labhist} have the underlying labeled topology depicted in Fig.~\ref{configa}A. The number of labeled histories of size $n$ with the same labeled topology $t$ is
\begin{equation}
\label{numblab}
\frac{(n-1)!}{\prod_{r=3}^n (r-1)^{d_r(t)}},
\end{equation}
where $d_r(t)$ is the number of internal nodes of $t$ from which exactly $r$ taxa descend \citep[][p.~46]{Steel16}. Eq.~(\ref{numblab}) also appears as the so-called ``shape functional'' of binary search trees \citep{Fill96}.

By summing the probability $1/|H_n|$ of each uniformly distributed labeled history of size $n$ with a given underlying labeled topology, the uniform distribution over the set $H_n$ induces the Yule--Harding (or Yule) distribution over the set $T_n$ of labeled topologies \citep{Yule25, Harding71, Brown94, McKenzieAndSteel00, SteelAndMcKenzie01, Rosenberg06:anncomb, ChangAndFuchs10, DisantoEtAl13, DisantoAndWiehe13}. The probability of a labeled topology $t$ can be calculated as
\begin{equation}
\label{eqPYH}
\text{P}_{\text{YH}}(t) = \frac{2^{n-1}}{n! \prod_{r=3}^n (r-1)^{d_r(t)} }.
\end{equation}
Under this distribution, among all labeled topologies with size $n$, those with the largest number of labeled histories have the highest probability. For balanced labeled topologies, the product in the denominator of Eq.~(\ref{eqPYH}) tends to be smaller than for unbalanced topologies, resulting in a greater probability.

% In Section \ref{yule}, by using the equivalence with uniformly distributed labeled histories mentioned above, we will derive distributional properties of the variable
% ``number of root configurations'' considered over labeled topologies of fixed size selected at random under the Yule--Harding distribution.

%%%%%%%%%%%%%%%%%%%%%%%%%%%%%%%%%%%%%%%%%%%%%%%%%%%%%%%%%%%%%%
%%%%%%%%%%%%%%%%%%%%%%%%%%%%%%%%%%%%%%%%%%%%%%%%%%%%%%%%%%%%%%
%%%%%%%%%%%%%%%%%%%%%%%%%%%%%%%%%%%%%%%%%%%%%%%%%%%%%%%%%%%%%%

\subsection{Asymptotic growth and analytic combinatorics}
\label{analyticomb}

Our study concerns the growth of increasing sequences. A sequence of non-negative numbers $a_n$ is said to have exponential growth $k^n$ or, equivalently, to be of exponential order $k$, if $a_n = k^n s(n)$, where $s$ is subexponential, that is, $\limsup_{n \rightarrow \infty} [s(n)^{1/n}] = 1$. Sequence $a_n$ grows exponentially in $n$ if its exponential order exceeds 1.

If $(a_n)$ has exponential order $k_a$ and $(b_n)$ has exponential order $k_b < k_a$, then the sequence of ratios $b_n/a_n$ converges to $0$ exponentially fast as $(k_b/k_a)^n$. If sequences $a_n$ and $b_n$ have the same exponential order, then we write $a_n \bowtie b_n$. If in addition the ratio $b_n/a_n$ converges to $1$, then we write $a_n \sim b_n$ and say that $(a_n)$ and $(b_n)$ have the same asymptotic growth.

Some results will make use of techniques of analytic combinatorics (see Sections IV and VI of \cite{FlajoletAndSedgewick09}). In particular, the entries of a sequence of integers $(a_n)_{n \geq 0}$ can be interpreted as coefficients of the power series expansion $A(z) = \sum_{n=0}^{\infty} a_n z^n$ at $z=0$ of a function $A(z)$, the generating function of the sequence. Considering $z$ as a complex variable, the behavior of $A(z)$ near its singularities---the points in the complex plane where $A(z)$ is not analytic---can provide information on the growth of its coefficients. Under suitable conditions, a correspondence exists between the expansion $A_{\alpha}(z)$ of the generating function $A(z)$ near its dominant singularity $\alpha$---that is, the singularity of smallest modulus---and the asymptotic growth of the coefficients $a_n$. In the simplest case, if $\alpha$ is the only dominant singularity of $A(z)$, then the $n$th coefficient $a_n$ of $A(z)$ has asymptotic growth $[z^n] A_{\alpha}(z)$, that is, the $n$th coefficient of $A_{\alpha}(z)$ (Theorem VI.4 of \cite{FlajoletAndSedgewick09}). In symbols,
$$a_n \sim [z^n] A_{\alpha}(z).$$
The exponential order of sequence $(a_n)$ is the inverse of the modulus of the dominant singularity $\alpha$ of $A(z)$ (Theorem IV.7 of \cite{FlajoletAndSedgewick09}). That is,
$$a_n \bowtie \alpha^{-n}.$$
As an example, sequence $|T_n|/n!$, with $|T_n|$ as in Eq.~(\ref{carciofo}), has exponential order 2 because $\alpha = \frac{1}{2}$ is the dominant singularity of the associated generating function in Eq.~(\ref{expo}). Thus, as $n \rightarrow \infty$, $|T_n|/n!$ increases with a subexponential multiple of $2^n$.

%%%%%%%%%%%%%%%%%%%%%%%%%%%%%%%%%%%%%%%%%%%%%%%%%%%%%%%%%%%%%%
%%%%%%%%%%%%%%%%%%%%%%%%%%%%%%%%%%%%%%%%%%%%%%%%%%%%%%%%%%%%%%
%%%%%%%%%%%%%%%%%%%%%%%%%%%%%%%%%%%%%%%%%%%%%%%%%%%%%%%%%%%%%%

\subsection{Ancestral configurations for matching gene trees and species trees}
\label{igno}

In this section, following \cite{DisantoAndRosenberg17}, we review features of the objects on which our study focuses: the ancestral configurations of a gene tree $G$ in a species tree $S$. In our framework, exactly one gene lineage has been selected from each species, and we assume $G$ and $S$ have the same labeled topology $t$.

\subsubsection{Definition of ancestral configurations}

Suppose $R$ is a realization of a gene tree $G$ in a species tree $S$, where $G=S=t$ (Fig.~\ref{configa}); $R$ is one of the possibilities for evolution of gene tree $G$ on matching species tree $S$. Looking backward in time, for node $\eta$ of $S$, consider the set $C(\eta,R)$ of gene lineages---edges of $G$---that are present in $S$ at the point just before node $\eta$.

The set $C(\eta,R)$ is the \emph{ancestral configuration} of $G$ at node $\eta$ of $S$. For example, for tree $t$ in Fig.~\ref{configa}A, with the realization $R_1$ of gene tree $G=t$ in the species tree $S=t$ in Fig.~\ref{configa}B, just before the root node $k$, the gene lineages present in the species tree are lineages $g$, $h$, and $i$. Hence, at species tree node $k$, the ancestral configuration is the set of gene lineages $C(k,R_1) = \{ g,h,i \}$. Similarly, the ancestral configuration of the gene tree at species tree node $j$ is $C(j,R_1)=\{ g,c,d \}$. In Fig.~\ref{configa}C, with a different realization $R_2$ of the same gene tree, the ancestral configuration at the species tree root $k$ is $C(k,R_2)=\{j,e,f \}$. The ancestral configuration at node $j$ is $C(j,R_2) = \{a,b,h \}$.

Let $\Re(G,S)$ be the set of realizations of gene tree $G=t$ in species tree $S=t$. For a given node $\eta$ of $t$, considering all possible elements $R \in \Re(G,S)$, the set of ancestral configurations is
\begin{equation}
\label{ciccio}
C(\eta) = \{ C(\eta,R) : R \in \Re(G,S) \}.
\end{equation}
The associated number of ancestral configurations is
\begin{equation}
\label{cicci}
c_{\eta} = |C(\eta)|.
\end{equation}
The quantity $c_{\eta}$ counts the ways the lineages of $G$ can reach the point right before node $\eta$ in $S$, considering all possible realizations of gene tree $G$ in species tree $S$. Choosing $t$ as in Fig.~\ref{configa}A, we have $C(g) = \{\{ a,b \}\}$, $C(h) = \{\{ c,d \}\}$, $C(i) = \{\{ e,f \}\}$, $C(j)=\{ \{a,b,c,d \}, \{g,c,d \}, \{a,b,h \},\{g,h \} \}$, and
\begin{equation}
\label{cm}
\small{
C(k) = \{ \{j,i \}, \{j,e,f \}, \{g,h,i \}, \{g,h,e,f \}, \{a,b,h,i \}, \{a,b,h,e,f \}, \{g,c,d,i \}, \{g,c,d,e,f \} , \{a,b,c,d,i \}, \{a,b,c,d,e,f \} \}.
}
\end{equation}
For different realizations $R_1,R_2 \in \Re(G,S)$ and an internal node $\eta$, it need not be true that $C(\eta,R_1) \neq C(\eta,R_2)$.

We say that a leaf or a 1-taxon tree has no ancestral configurations. In addition, the definition of an ancestral configuration at node $\eta$, by considering the point right before node $\eta$ in the species tree, excludes the case in which all gene tree lineages descended from gene tree node $\eta$ have coalesced at species tree node $\eta$. Thus, $\{\eta \} \notin C(\eta)$.

Because we consider the case of $G=S=t$, the set $C(\eta)$ and the quantity $c_{\eta}$ in Eqs.~(\ref{ciccio}) and (\ref{cicci}) depend only on node $\eta$ and tree $t$. We use the term \emph{configurations at node $\eta$ of $t$} to denote elements of $C(\eta)$.

%%%%%%%%%%%%%%%%%%%%%%%%%%%%%%%%%%%%%%%%%%%%%%%%%%%%%%%%%%%%%%
%%%%%%%%%%%%%%%%%%%%%%%%%%%%%%%%%%%%%%%%%%%%%%%%%%%%%%%%%%%%%%
%%%%%%%%%%%%%%%%%%%%%%%%%%%%%%%%%%%%%%%%%%%%%%%%%%%%%%%%%%%%%%

\subsubsection{Root and total configurations}
\label{ignores}

Our focus is on configurations at the root of $t$. Let $N(t)$ be the set of nodes of a tree $t$, including both leaf nodes and internal nodes. With $|t|$ leaf nodes and $|t|-1$ internal nodes in $t$, $|N(t)|=2|t|-1$. Define the \emph{total} number of configurations in $t$ by
$$c = \sum_{\eta \in N(t)} c_{\eta}.$$
Let $c_r$ be the number of configurations at the root $r$ of $t$, or \emph{root configurations} for short. Because $c_r \geq c_{\eta}$ for each node $\eta$ of $t$, we have
\begin{equation}
\label{bau}
c_r \leq c \leq (2|t|-1) c_r.
\end{equation}
Quantities $c$ and $c_r$ are equal up to a factor that is at most polynomial in $|t|$, and they have the same exponential order when measured across families of trees of increasing size.

Selecting a tree of size $n$ at random from the set of labeled topologies, inequality (\ref{bau}) gives $\mathbb{E}_n[c_r] \leq \mathbb{E}_n[c] \leq 2n\mathbb{E}_n[c_r]$ and $\mathbb{E}_n[c_r^2] \leq \mathbb{E}_n[c^2] \leq 4n^2\mathbb{E}_n[c_r^2]$. In expectation $\mathbb{E}$ and variance $\mathbb{V}$, exponential growth for total configurations follows that for root configurations:
\begin{eqnarray} \nonumber
\mathbb{E}_n[c] & \bowtie & \mathbb{E}_n[c_r] \\ \nonumber
\mathbb{E}_n[c^2] & \bowtie & \mathbb{E}_n[c_r^2] \\ \nonumber
\mathbb{V}_n[c] = \mathbb{E}_n[c^2] - \mathbb{E}_n[c]^2 & \bowtie & \mathbb{E}_n[c_r^2] - \mathbb{E}_n[c_r]^2 = \mathbb{V}_n[c_r].
\end{eqnarray}

%%%%%%%%%%%%%%%%%%%%%%%%%%%%%%%%%%%%%%%%%%%%%%%%%%%%%%%%%%%%%%
%%%%%%%%%%%%%%%%%%%%%%%%%%%%%%%%%%%%%%%%%%%%%%%%%%%%%%%%%%%%%%
%%%%%%%%%%%%%%%%%%%%%%%%%%%%%%%%%%%%%%%%%%%%%%%%%%%%%%%%%%%%%%

\subsubsection{Known results}
\label{ignoreste}

We recall some results of \cite{DisantoAndRosenberg17} on the number of configurations possessed by a tree.

\smallskip

(i) For a given tree $t$ with $|t|>1$, let $r$ denote the root node of $t$, with $r_{L}$ and $r_R$ being the two child nodes of $r$. The number $c_r$ of possible configurations at $r$ can be recursively computed as
\begin{equation}
\label{eqC}
c_r = (c_{r_L}+1)(c_{r_R}+1),
\end{equation}
where we set $c_r = 0$ if $|t|=1$. For example, for the tree of Fig.~\ref{configa}A, we have $r=k$, $r_L=j$, $r_R=i$, and $c_k=10=(4+1)(1+1)=(c_j+1)(c_i+1)$, as determined by Eq.~(\ref{eqC}).

\smallskip

(ii) Consider a representative labeling of each unlabeled topology of size $n$. Among these trees, the largest number of root configurations and the largest total number of configurations have exponential order $k_0$, where $k_0 \approx 1.5028$. The smallest number of root configurations and the smallest total number of configurations have polynomial growth with the tree size $n$. Furthermore, consider the balanced family of unlabeled topologies defined recursively by $|t_1|=1$ and $t_n = (t_d, t_{n-d})$, where $d$ denotes the power of $2$ nearest to $\frac{n}{2}$. Among the unlabeled topologies with $n$ taxa, $t_n$ has the largest number of root configurations. The maximally asymmetric caterpillar unlabeled topology has the smallest number of root configurations.

\smallskip

(iii) For a labeled topology of given size $n$ selected uniformly at random, the mean number of root configurations $c_r$ and the mean total number of configurations $c$ grow asymptotically like
\begin{eqnarray}
\label{meandis}
\mathbb{E}_n[c_r] & \sim    & \sqrt{\frac{3}{2}} \left(\frac{4}{3}\right)^n, \\ \label{meandiss}
\mathbb{E}_n[c]   & \bowtie & \left(\frac{4}{3}\right)^n.
\end{eqnarray}
The variances of $c_r$ and $c$ satisfy the asymptotic relations
\begin{eqnarray}
\label{vardis}
\mathbb{V}_n[c_r] & \sim    & \sqrt{\frac{7(11 - \sqrt{2})}{34}}\left[\frac{4}{7 (8\sqrt{2} - 11) }\right]^n, \\ \label{vardiss}
\mathbb{V}_n[c]   & \bowtie & \left[\frac{4}{7 (8\sqrt{2} - 11) }\right]^n .
\end{eqnarray}

%%%%%%%%%%%%%%%%%%%%%%%%%%%%%%%%%%%%%%%%%%%%%%%%%%%%%%%%%%%%%%
%%%%%%%%%%%%%%%%%%%%%%%%%%%%%%%%%%%%%%%%%%%%%%%%%%%%%%%%%%%%%%
%%%%%%%%%%%%%%%%%%%%%%%%%%%%%%%%%%%%%%%%%%%%%%%%%%%%%%%%%%%%%%

\subsection{Additive tree parameters and root configurations}
\label{addp}

A quantity $F(t)$ that is computed for a tree $t$ and whose value can be calculated as
\begin{equation*}
%\label{addt}
F(t) = F(t_L) + F(t_R) + f(t),
\end{equation*}
where $t_L$ and $t_R$ are the two root subtrees of $t$, is called an \emph{additive tree parameter} with \emph{toll function} $f(t)$ \citep[e.g.][]{Wagner15}. For a variety of tree families, \cite{Wagner15} showed that an additive tree parameter $F(t)$ is asymptotically normally distributed if the toll function $f(t)$ is bounded and the mean value of $|f(t)|$, considered over uniformly distributed trees of fixed size, goes to $0$ exponentially fast as the tree size increases.

For a tree $t$, consider the quantity $\log(c_r + 1)$, that is, the natural logarithm of one more than the number of root configurations of $t$. From Eq.~(\ref{eqC}), a simple calculation yields for $|t| \geq 2$
\begin{equation}
\label{eqlogC}
\log(c_r + 1) = \log(c_{r_L}+1) + \log(c_{r_R}+1) + \log\left( 1+\frac{1}{c_r} \right).
\end{equation}
In Eq.~(\ref{eqlogC}), if we set
$$F(t) = \log[c_r(t) + 1],$$
then the associated toll function is given for $|t| \geq 2$ by
$$f(t) = \log\left[ 1+\frac{1}{c_r(t)} \right].$$
We set $f(t) = F(t) = \log(1) = 0$ if $|t| = 1$. We can therefore consider root configurations in the context of additive tree parameters.

% Using results of \cite{Wagner15} for the families of ordered unlabeled topologies (pruned binary trees in \cite{Wagner15}) and ordered unlabeled histories (binary
% increasing trees in \cite{Wagner15}), in Sections~\ref{unif} and \ref{yule}, we show that under the uniform and Yule--Harding distributions over labeled topologies of size % $n$, the random variable ``logarithm of the number of root configurations,'' suitably rescaled, converges to a normal distribution. In Section~\ref{yule}, we also study
% the mean and variance of the number of root configurations for labeled topologies of fixed size under the Yule--Harding distribution.

%%%%%%%%%%%%%%%%%%%%%%%%%%%%%%%%%%%%%%%%%%%%%%%%%%%%%%%%%%%%%%
%%%%%%%%%%%%%%%%%%%%%%%%%%%%%%%%%%%%%%%%%%%%%%%%%%%%%%%%%%%%%%
%%%%%%%%%%%%%%%%%%%%%%%%%%%%%%%%%%%%%%%%%%%%%%%%%%%%%%%%%%%%%%

\section{Equivalences for the distribution of the number of root configurations}
\label{secEquivalences}

We prove a series of equivalences needed for analyzing distributional properties of the number of root configurations. In Section \ref{pinto}, we show that the distribution of the number of root configurations over uniformly distributed labeled topologies or labeled histories can be analyzed by considering equivalently the distribution of the number of root configurations over uniformly distributed ordered unlabeled topologies or ordered unlabeled histories, respectively. In Section \ref{pruned}, we obtain a correspondence between antichains of pruned binary trees and root configurations of ordered unlabeled topologies.

%%%%%%%%%%%%%%%%%%%%%%%%%%%%%%%%%%%%%%%%%%%%%%%%%%%%%%%%%%%%%%
%%%%%%%%%%%%%%%%%%%%%%%%%%%%%%%%%%%%%%%%%%%%%%%%%%%%%%%%%%%%%%
%%%%%%%%%%%%%%%%%%%%%%%%%%%%%%%%%%%%%%%%%%%%%%%%%%%%%%%%%%%%%%

\subsection{Equivalences with ordered unlabeled topologies and histories}
\label{pinto}

Distributional properties of a tree parameter defined over the set of labeled topologies can in some cases be investigated by studying the same parameter over a different tree family. In particular, if the tree parameter under consideration depends only on tree topology, then its distribution can be equivalently analyzed over a different tree set taken under a probability model that induces or is induced by the probability model assumed for labeled topologies. In this direction, \cite{BlumEtAl06} derived a general framework for analyzing tree parameters of labeled topologies under a variety of probabilistic models defined over binary search trees.

In this section, we obtain results analogous to those of \cite{BlumEtAl06}. We show that the number of root configurations---or any other tree parameter that depends only on the branching structure of the tree---has the same distribution when considered over uniformly distributed labeled topologies or over uniformly distributed ordered unlabeled topologies of the same size (Lemma \ref{lem}). Similarly, the number of root configurations has the same distribution over uniformly distributed labeled histories of size $n$ as for uniformly distributed ordered unlabeled histories of size $n$ (Lemma \ref{lema}).

Moreover, because the uniform distribution over the set of labeled histories of size $n$ induces the Yule--Harding distribution over the set of labeled topologies of size $n$ (Section \ref{yhdistr}), as a direct consequence of Lemma \ref{lema} we have that the number of root configurations has the same distribution when considered over Yule--Harding-distributed labeled topologies or over uniformly distributed ordered unlabeled histories (Lemma \ref{corr}). By using these facts, Propositions \ref{leeemRn} and \ref{lemRn} give recursive formulas for the probabilities under the uniform and Yule--Harding probability models, respectively, that a random labeled topology of size $n$ has $c_r=\rho$ root configurations.

\begin{lemm}
\label{lem}
The distribution of the number of root configurations over labeled topologies of size $n$ selected uniformly at random matches the distribution of the number of root configurations over ordered unlabeled topologies of size $n$ selected uniformly at random.
\end{lemm}

\noindent \emph{Proof.} First, we note that the number of root configurations of a labeled topology or ordered unlabeled topology depends only on the underlying unlabeled topology. Thus, to prove the claim, it suffices to show that for each unlabeled topology $t$ of size $n$, we have
\begin{equation}
\label{eqOrlab}
\frac{\text{or}(t)}{C_{n-1}} = \frac{\text{lab}(t)}{|T_n|},
\end{equation}
where $\text{or}(t)$ and $\text{lab}(t)$ are the number of orientations of $t$ and the number of leaf labelings of $t$, respectively. Note from Eqs.~(\ref{cata}) and (\ref{carciofo}) that ${\text{or}(t)}/{C_{n-1}}$ and $ {\text{lab}(t)}/{|T_n|}$ give the probability of the unlabeled topology $t$ induced by the uniform distribution over the set of ordered unlabeled topologies and labeled topologies of $n$ taxa, respectively.

By using $C_{n-1} = {{2n-2}\choose{n-1}}/n $ and $|T_n| = (2n-2)!/[2^{n-1} (n-1)!] $ from Eqs.~(\ref{cata}) and (\ref{carciofo}), Eq.~(\ref{eqOrlab}) can be rewritten
\begin{equation*}
\label{eqlabt}
\text{lab}(t) = \text{or}(t) \frac{n!}{2^{n-1}},
\end{equation*}
which we demonstrate by induction on the size of $t$. Let $t_L$ and $t_R$ be the two root subtrees of $t$, with sizes $|t_L| = L$ and $|t_R|=R$. Thus, for $n \geq 2$,
\begin{eqnarray}
\label{lem1}
\text{lab}(t) & = & \text{lab}(t_L) \, \text{lab}(t_R) \, {{n}\choose{L}} \, \frac{1}{1+ \delta_{t_L=t_R}} \\
\text{or}(t)  & = & \text{or}(t_L)  \, \text{or}(t_R)  \, \frac{2}{1+ \delta_{t_L=t_R}},
\end{eqnarray}
where $\delta_{t_L=t_R}=1$ if $t_L = t_R$, and $\delta_{t_L=t_R}=0$ otherwise. If we insert $ \text{lab}(t_L) = \text{or}(t_L) {L!}/{2^{L-1}}$ and $\text{lab}(t_R) = \text{or}(t_R) {R!}/{2^{R-1}}$ into Eq.~(\ref{lem1}), then we find
\begin{eqnarray}
\text{lab}(t) &=& \text{or}(t_L) \, \text{or}(t_R) \, \frac{L! \, R!}{2^{n-2}} {{n}\choose{L}} \, \frac{1}{1+ \delta_{t_L=t_R}} \\
&=& \text{or}(t_L) \, \text{or}(t_R) \, \frac{n!}{2^{n-1}} \, \frac{2}{1+ \delta_{t_L=t_R}} = \text{or}(t) \frac{n!}{2^{n-1}},
\end{eqnarray}
as desired.
$\Box$

\medskip

The proof shows that the ratio of orderings to labelings for an unlabeled topology is independent of the unlabeled topology. Hence, because the number of root configurations of a labeled topology or ordered unlabeled topology depends only on the underlying unlabeled topology, the probability that a labeled topology chosen uniformly at random has $\rho$ root configurations equals the probability that an ordered unlabeled topology chosen uniformly at random has $\rho$ root configurations. We use Lemma \ref{lem} to calculate the probability that a labeled topology of size $n$ selected under the uniform distribution has $\rho$ root configurations as the probability that an ordered unlabeled topology of size $n$ selected under the uniform distribution has $\rho$ root configurations.

\begin{prop}
\label{leeemRn}
Let $R_n$ be the random variable that represents the number of root configurations in an ordered unlabeled topology of size $n$ selected uniformly at random. (i) We have $R_1=0$, and for $n \geq 2$,
\begin{equation}
\label{rectot}
R_n \stackrel{d}{=} (R_{I_n}+1)( R^{*}_{n-I_n}+1),
\end{equation}
where $I_n$ is distributed over the interval $[1,n-1]$ with Catalan probability $\mathbb{P}[I_n=j]={C_{j-1}C_{n-j-1}}/{C_{n-1}}$, $R_j^{*}$ is an independent copy of $R_j$ for each $j \in [1,n-1]$, and both $R_j$ and $R^{*}_j$ are independent of $I_j$ for $j \in [1,n-1]$. Furthermore, (ii) the probability that a random labeled topology of size $n$ selected under the uniform distribution has $c_r = \rho$ root configurations can be calculated as $\mathbb{P}[c_r = \rho] = \mathbb{P}[R_n = \rho]$, where $\mathbb{P}[R_n = \rho]$ has recursive formula
\begin{equation}
\label{esattaprobuni}
\mathbb{P}[R_n = \rho] = \sum_{d \in \text{Div}(\rho)} \sum_{j=1}^{n-1} \mathbb{P}[I_n = j] \, \mathbb{P}[R_j=d-1] \, \mathbb{P}\bigg[ R_{n-j} = \frac{\rho}{d} - 1 \bigg],
\end{equation}
$\text{Div}(\rho)$ denotes the set of positive integers that divide $\rho$, $\mathbb{P}[I_n = j] = {C_{j-1} C_{n-j-1}}/{C_{n-1}}$, and $\mathbb{P}[R_n = 0] = \delta_{n,1}$.
\end{prop}

\noindent \emph{Proof.}
The recurrence in Eq.~(\ref{rectot}) follows from Eq.~(\ref{eqC}). Observe that for a random uniform ordered unlabeled topology $t$ of $n$ taxa, the probability that the left (or right) root subtree of $t$ has size $I_n = j$ is given by $\mathbb{P}[I_n=j] = {C_{j-1} C_{n-j-1}}/{C_{n-1}}$, where $C_{j-1}$, $C_{n-j-1}$, and $C_{n-1}$ give the numbers of ordered unlabeled topologies of size $j$, $n-j$, and $n$, respectively (Section \ref{out}). This establishes (i).

For (ii), Eq.~(\ref{esattaprobuni}) is a direct consequence of Lemma \ref{lem} and Eq.~(\ref{rectot}).
$\Box$

\medskip

We now consider the equivalence between uniformly distributed labeled histories and uniformly distributed ordered unlabeled histories.

\begin{lemm}
\label{lema}
The distribution of the number of root configurations over labeled histories of size $n$ selected uniformly at random matches the distribution of the number of root configurations over ordered unlabeled histories of size $n$ selected uniformly at random.
\end{lemm}

\noindent \emph{Proof.} The proof is similar to that of Lemma \ref{lem}: we show that for each unlabeled history $t$ of size $n$, we have
\begin{equation}
\label{eqFH}
\frac{\text{or}(t)}{F_{n-1}} = \frac{\text{lab}(t)}{ |H_n| },
\end{equation}
where $\text{or}(t)$ and $\text{lab}(t)$ are the number of orientations of $t$ and the number of leaf labelings of $t$, respectively. In other words, we prove that the uniform distribution over the set of ordered unlabeled histories of size $n$ and the uniform distribution over the set of labeled histories of size $n$ both induce the same probability distribution over the set of unlabeled histories of $n$ taxa. The same property has already been shown by Lambert and Stadler (2013, p.~116), \nocite{LambertAndStadler13} following a slightly different approach.

Using $F_{n-1} = (n-1)!$ and $|H_n| = n! (n-1)! / 2^{n-1}$ from Eqs.~(\ref{fact}) and (\ref{carciofino}), Eq.~(\ref{eqFH}) can be rewritten
$$\text{lab}(t) = \text{or}(t) \frac{n!}{2^{n-1}},$$
which we verify by induction on $|t|$. Let $t_L$ and $t_R$ denote the two root subtrees of $t$, with sizes $|t_L|=L$ and $|t_R|=R$. Hence, for $n \geq 2$ we have
\begin{eqnarray}
\label{lem2}
\text{lab}(t) & = & \text{lab}(t_L) \, \text{lab}(t_R) \, {{n}\choose{L}} \\
\text{or}(t)  & = & 2 \, \text{or}(t_L) \, \text{or}(t_R).
\end{eqnarray}
By setting $ \text{lab}(t_L) = \text{or}(t_L) {L!}/{2^{L-1}}$ and $\text{lab}(t_R) = \text{or}(t_R) {R!}/{2^{R-1}}$ in Eq.~(\ref{lem2}), we find
\begin{eqnarray}
\text{lab}(t) &=& \text{or}(t_L) \, \text{or}(t_R) \, \frac{L! \, R!}{2^{n-2}} {{n}\choose{L}} \\
              &=& \text{or}(t_L) \, \text{or}(t_R) \, \frac{2 \, n!}{2^{n-1}}  =  \text{or}(t) \frac{n!}{2^{n-1}},
\end{eqnarray}
as desired.
$\Box$

\medskip

Next, we translate the result of Lemma \ref{lema} in terms of Yule--Harding-distributed labeled topologies.

\begin{lemm}
\label{corr}
The distribution of the number of root configurations over labeled topologies of size $n$ selected according to the Yule--Harding distribution matches the distribution of the number of root configurations over ordered unlabeled histories of size $n$ selected uniformly at random.
\end{lemm}

\noindent \emph{Proof.} The equivalence follows from Lemma \ref{lema} and the fact that the uniform distribution over labeled histories of size $n$ induces the Yule--Harding distribution on the set of labeled topologies of size $n$ (Section \ref{yhdistr}). $\Box$

\medskip

By Lemma \ref{corr}, we can calculate the probability that a labeled topology of size $n$ selected under the Yule--Harding distribution has $\rho$ root configurations as the probability that a random uniform ordered unlabeled history of size $n$ has $\rho$ root configurations. In particular, we have the following proposition.

\begin{prop}
\label{lemRn}
Let $R_n$ be the random variable that represents the number of root configurations in an ordered unlabeled history of size $n$ selected uniformly at random. (i) We have $R_1 =0$, and for $n \geq 2$,
\begin{equation}
\label{eqRn}
R_n \stackrel{d}{=} (R_{I_n}+1)( R^{*}_{n-I_n}+1),
%= R_{I_n} \, R^{*}_{n-I_n} + R_{I_n} + R^{*}_{n-I_n} + 1,
\end{equation}
where $I_n$ is uniformly distributed over the interval $[1,n-1]$, $R^{*}_{j}$ is an independent copy of $R_j$ for each $j \in [1,n-1]$, and both $R_j$ and $R^{*}_j$ are independent of $I_j$ for $j \in [1,n-1]$. Furthermore, (ii) the probability that a random labeled topology of size $n$ selected under the Yule--Harding distribution has $c_r = \rho$ root configurations can be calculated as $\mathbb{P}[c_r = \rho] = \mathbb{P}[R_n = \rho]$, where $\mathbb{P}[R_n = \rho]$ has recursive formula
\begin{equation}
\label{esattaprobyul}
\mathbb{P}[R_n = \rho] = \sum_{d \in \text{Div}(\rho)} \sum_{j=1}^{n-1} \mathbb{P}[I_n = j] \, \mathbb{P}[R_j=d-1] \, \mathbb{P}\bigg[R_{n-j}= \frac{\rho}{d}-1 \bigg],
\end{equation}
$\text{Div}(\rho)$ denotes the set of positive integers that divide $\rho$, $\mathbb{P}[I_n=j] = \frac{1}{n-1}$, and $\mathbb{P}[R_n=0] = \delta_{n,1}$.
\end{prop}
\noindent \emph{Proof.} The formula in Eq.~(\ref{eqRn}) follows directly from Eq.~(\ref{eqC}) when we observe that, for a random uniform ordered unlabeled history $t$ of $n$ taxa, the probability that the left (or right) root subtree of $t$ has size $I_n = j$ is
$$\mathbb{P}[I_n=j] = \frac{F_{j-1} F_{n-j-1} {{n-2}\choose{j-1}}}{F_{n-1}} = \frac{1}{n-1}.$$
Eq.~(\ref{esattaprobyul}) is a direct consequence of Lemma \ref{corr} and Eq.~(\ref{eqRn}). $\Box$

%%%%%%%%%%%%%%%%%%%%%%%%%%%%%%%%%%%%%%%%%%%%%%%%%%%%%%%%%%%%%%
%%%%%%%%%%%%%%%%%%%%%%%%%%%%%%%%%%%%%%%%%%%%%%%%%%%%%%%%%%%%%%
%%%%%%%%%%%%%%%%%%%%%%%%%%%%%%%%%%%%%%%%%%%%%%%%%%%%%%%%%%%%%%

\subsection{Equivalences with antichains of pruned binary trees}
\label{pruned}

To use results of \cite{Wagner15} to obtain probability distributions for root configurations, we must translate between root configurations for labeled topologies and non-empty antichains for pruned binary trees.

A pruned binary tree is an ordered unlabeled topology in which the external branches---those terminating in a leaf---have been removed. To illustrate the pruning operation, consider the ordered unlabeled topology depicted on the left of Fig.~\ref{orient}A and assign arbitrary labels to all its nodes, as in Fig.~\ref{configa}A. The leaf labels of the pruned binary tree resulting from this process can be described by the Newick format $((g,h),i)$. Note that pruned binary trees have their left--right orientation induced by the overlying ordered unlabeled topology.

If $t$ is an ordered unlabeled topology of size $n$ and $\tilde{t}$ is its associated pruned binary tree of $n-1$ nodes, then we can consider $\tilde{t}$ as the Hasse diagram of a partially ordered set with ground set given by the nodes of $\tilde{t}$---the \emph{internal} nodes of $t$---and order relation determined by the descendant--ancestor relationship in $\tilde{t}$. An antichain of $\tilde{t}$ is a subset of its nodes such that no two elements in the subset are comparable by the order relation. For instance, the two-element antichains of pruned binary tree $((g,h),i)$ in Fig.~\ref{configa}A are $\{g,h\}, \{g,i\}, \{h,i\},$ and $\{j,i\}$.

The non-empty antichains of the pruned binary tree $\tilde{t}$ bijectively correspond to the root configurations of the overlying ordered unlabeled topology $t$: omitting leaves from a root configuration of $t$ yields an antichain of $\tilde{t}$, and adding leaves to an antichain of $\tilde{t}$ so that each leaf of $t$ is either represented or has one of its ancestral nodes represented yields a root configuration of $t$.

For instance, consider the set in Eq.~(\ref{cm}) of the root configurations of the ordered unlabeled topology in Fig.~\ref{configa}A. By omitting leaves from each configuration, we obtain the antichains of $\tilde{t}$:
\begin{equation*}
%\label{cmant}
\small{
\{ \{j,i \}, \{j \}, \{g,h,i \}, \{g,h \}, \{h,i \}, \{h \}, \{g,i \}, \{g \} , \{i \}, \emptyset \}.
}
\end{equation*}
We make a substitution of the empty antichain $\emptyset$ that emerges from the root configuration consisting of all the leaves by the antichain $\{ k \}$ consisting only of the root of $\tilde{t}$; we have then bijectively paired all root configurations of $t$ and all non-empty antichains of $\tilde{t}$. Using this correspondence, we have the next result.
\begin{lemm}
\label{anticha}
The distribution of the number of root configurations over labeled topologies of size $n$ selected uniformly at random matches the distribution of the number of non-empty antichains over the set of $(n-1)$-node pruned binary trees selected uniformly at random.
\end{lemm}

\noindent \emph{Proof.} By Lemma \ref{lem}, the number of root configurations has the same distribution when considered over uniformly distributed labeled topologies of size $n$ or over uniformly distributed ordered unlabeled topologies of size $n$. By the correspondence between antichains of pruned binary trees with $n-1$ nodes and root configurations of associated ordered unlabeled topologies of size $n$, the distribution of the number of root configurations over uniformly distributed ordered unlabeled topologies of size $n$ matches the distribution of the number of non-empty antichains over uniformly distributed pruned binary trees with $n-1$ nodes. $\Box$

%%%%%%%%%%%%%%%%%%%%%%%%%%%%%%%%%%%%%%%%%%%%%%%%%%%%%%%%%%%%%%
%%%%%%%%%%%%%%%%%%%%%%%%%%%%%%%%%%%%%%%%%%%%%%%%%%%%%%%%%%%%%%
%%%%%%%%%%%%%%%%%%%%%%%%%%%%%%%%%%%%%%%%%%%%%%%%%%%%%%%%%%%%%%

\section{Root configurations under the uniform distribution on labeled topologies}
\label{unif}

\cite{DisantoAndRosenberg17} determined the mean and variance of the number of root configurations for uniformly distributed labeled topologies of size $n$ (Section \ref{ignoreste}). In this section, we use the correspondence with antichains given in Section \ref{pruned} to show that the logarithm of the number of root configurations for uniformly distributed labeled topologies of size $n$, suitably rescaled, converges to a normal distribution.

Wagner (2015, Section 2.3.2) \nocite{Wagner15} studied the number $a(t)$ of non-empty antichains of a randomly selected pruned binary tree $t$ of given size. For a pruned binary tree of $n$ nodes selected uniformly at random, he considered $\log a(t)$, showing that $\big(\log a - \mathbb{E}_n[\log a ]\big)/\sqrt{\mathbb{V}_n[\log a ]}$ converges to a standard normal distribution as $n \rightarrow \infty$, where $\mathbb{E}_n[\log a] \sim \mu n$ and $\mathbb{V}_n[\log a] \sim \sigma^2 n$, with constants $(\mu, \sigma^2) \approx (0.272, 0.034)$.

By Lemma \ref{anticha}, Wagner's variable $\log a$ asymptotically has the same distribution as the variable $\log c_r$ considered over uniformly distributed labeled topologies of size $n+1$. We thus have the following result.
\begin{prop}
\label{lognormunil}
The logarithm of the number of root configurations in a labeled topology of size $n$ selected uniformly at random, rescaled as $\big(\log c_r - \mathbb{E}_n[\log c_r]\big)/\sqrt{\mathbb{V}_n[\log c_r]}$, converges to a standard normal distribution, where $\mathbb{E}_n[\log c_r] \sim \mu n$ and $\mathbb{V}_n[\log c_r] \sim \sigma^2 n$, $(\mu, \sigma^2) \approx (0.272, 0.034)$.
\end{prop}

The result gives an asymptotic lognormal distribution for the number of root configurations of a labeled topology of size $n$ selected uniformly at random. Although we do not expect $e^{\mathbb{E}_n[\log c_r]}$ and $e^{\sigma_n[\log c_r]}$ to agree with $\mathbb{E}_n[c_r]$ and $\sigma_n[c_r]$, for the mean we see that in the $n \rightarrow \infty$ limit, $e^{\mathbb{E}_n[\log c_r]} \approx e^{0.272n} \approx 1.313^n$, numerically close to the exponential growth of $\mathbb{E}_n[c_r]$, or $(4/3)^n$ (Eq.~(\ref{meandis})). For, the standard deviation $e^{\sigma_n[\log c_r]} \approx e^{\sqrt{0.034}n} \approx 1.202^n$ is not as close to the exponential growth of $\sigma_n[c_r]$ from Eq.~(\ref{vardis}), which gives $[2/\sqrt{7(8\sqrt{2}-11)}]^n \approx 1.350^n$.

\begin{figure}
\begin{center}
\includegraphics*[scale=0.50,trim=0 0 0 0]{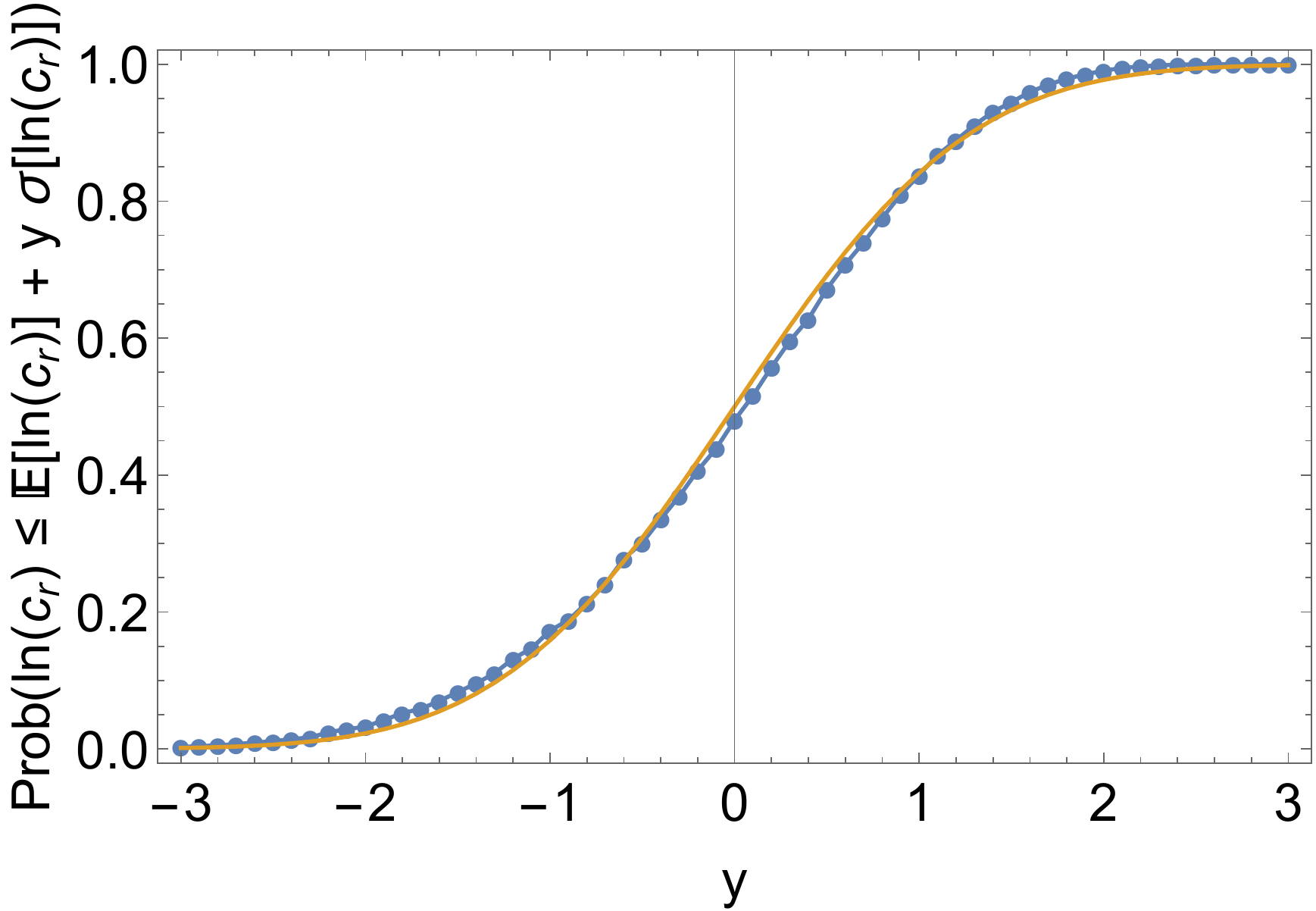}
\end{center}
\vspace{-.7cm}
\caption{{\small Cumulative distribution of the natural logarithm of the number of root configurations for uniformly distributed labeled topologies of size $n=15$ (dotted line). Each dot has its abscissa determined by a value of $y$ ranging in the interval $y \in [-3,3]$ in steps of $0.1$. Given $y$, the quantity plotted is the probability that a labeled topology with $n=15$ chosen uniformly at random has a number of root configurations less than or equal to $\exp\left( \mathbb{E}[\log c_r] + y \sigma[\log c_r] \right)$, where $\mathbb{E}[\log c_r]$ and $\sigma[\log c_r$] are respectively the mean and standard deviation of the logarithm of the number of root configurations for uniformly distributed labeled topologies with $n=15$ taxa (Proposition \ref{lognormunil}). The solid line is the cumulative distribution of a Gaussian random variable with mean 0 and variance 1.
}}
\label{cumuni}
\end{figure}

For fixed $n$, we can compute the exact distribution of $c_r$ and $\log c_r$ under a uniform distribution across labeled topologies of size $n$, as described in Proposition \ref{leeemRn}ii. Fig.~\ref{cumuni} shows the cumulative distribution $\mathbb{P}\big[\log c_r \leq \mathbb{E}[\log c_r] + y \sigma[\log c_r] \big]$ as a function of $y$, when labeled topologies are selected uniformly at random among the $2.13 \times 10^{14}$ labeled topologies with 15 leaves. To obtain the distribution, we can count root configurations for arbitrary labelings of each of the 4850 unlabeled topologies with 15 leaves, and then count labelings for each unlabeled topology~\cite[][p.~47]{Steel16}. Already for small tree size, the figure shows that the exact cumulative distribution is close to the cumulative distribution of a Gaussian random variable with mean 0 and variance 1.

%%%%%%%%%%%%%%%%%%%%%%%%%%%%%%%%%%%%%%%%%%%%%%%%%%%%%%%%%%%%%%
%%%%%%%%%%%%%%%%%%%%%%%% Figure 4 %%%%%%%%%%%%%%%%%%%%%%%%%%%%
%%%%%%%%%%%%%%%%%%%%%%%%%%%%%%%%%%%%%%%%%%%%%%%%%%%%%%%%%%%%%%
%%%%%%%%%%%%%%%%%%%%%%%%%%%%%%%%%%%%%%%%%%%%%%%%%%%%%%%%%%%%%%

%%%%%%%%%%%%%%%%%%%%%%%%%%%%%%%%%%%%%%%%%%%%%%%%%%%%%%%%%%%%%%
%%%%%%%%%%%%%%%%%%%%%%%%%%%%%%%%%%%%%%%%%%%%%%%%%%%%%%%%%%%%%%
%%%%%%%%%%%%%%%%%%%%%%%%%%%%%%%%%%%%%%%%%%%%%%%%%%%%%%%%%%%%%%

\section{Root configurations under the Yule--Harding distribution on labeled topologies}
\label{yule}

We next study distributional properties of the number of root configurations for labeled topologies selected under the Yule--Harding probability model. Section \ref{yhdistr} noted that this model assigns higher probability to trees with a high degree of balance compared to that assigned by the uniform model; Section \ref{ignoreste} noted that balanced trees have high numbers of root configurations relative to unbalanced trees. We therefore find that the mean number of root configurations for labeled topologies of size $n$ grows exponentially faster under the Yule--Harding model than under the uniform model. The variance of the number of root configurations also has faster growth.

%%%%%%%%%%%%%%%%%%%%%%%%%%%%%%%%%%%%%%%%%%%%%%%%%%%%%%%%%%%%%%
%%%%%%%%%%%%%%%%%%%%%%%%%%%%%%%%%%%%%%%%%%%%%%%%%%%%%%%%%%%%%%
%%%%%%%%%%%%%%%%%%%%%%%%%%%%%%%%%%%%%%%%%%%%%%%%%%%%%%%%%%%%%%

\subsection{Lognormal distribution of the number of root configurations}
\label{yulelognormal}

We begin the analysis of the number of root configurations under the Yule--Harding distribution by showing that the logarithm of the number of root configurations of a Yule--Harding random labeled topology of size $n$, when suitably rescaled, converges to a standard normal distribution.

The results in this section are obtained by considering root configurations over ordered unlabeled histories of given size selected under the uniform distribution. Owing to Lemma \ref{corr}, we can demonstrate that the number of root configurations in a Yule--Harding random labeled topology of size $n$ asymptotically follows a lognormal distribution by showing that the number of root configurations is asymptotically lognormally distributed when considered over the set of uniformly distributed ordered unlabeled histories of $n$ taxa. We use a result of \cite{Wagner15} for additive tree parameters of ordered unlabeled histories. We first must verify a technical condition for the mean of the random variable $\log\left( 1 + 1/c_r \right)$, considered over uniformly distributed ordered unlabeled histories.
\begin{lemm}
\label{le}
For uniformly distributed ordered unlabeled histories of size $n$, the mean value $\mathbb{E}_n\left[ \log \left(1 + 1/c_r \right) \right]$ of the random variable $\log(1 + 1/c_r)$ converges to $0$ exponentially fast as $n$ increases. In particular,
\begin{equation}
\label{teclem}
\mathbb{E}_n\left[ \log \left(1 + \frac{1}{c_r} \right) \right] = \mathcal{O}(0.9^n).
\end{equation}
\end{lemm}

\noindent \emph{Proof.} To show that $\mathbb{E}_n\left[ \log \left(1 + 1/c_r \right) \right]$ has exponential growth $\mathcal{O}(0.9^n)$ for an ordered unlabeled history $t$ of size $n$ selected uniformly at random, we consider the mean value $\mathbb{E}_n[2^{-\text{ch}}]$ of the random variable $2^{-\text{ch}}$---where $\text{ch}$ is the number of cherries in $t$. We claim that
\begin{equation}
\label{cherr}
\mathbb{E}_n[2^{-\text{ch}}] = \mathcal{O}(0.9^n).
\end{equation}

For a tree $t$ with $|t| \geq 3$, $c_r(t) \geq 2^{\text{ch}(t)}$, as each cherry node generates a pair of ancestral configurations: the configuration corresponding to the node, and the configuration corresponding to its pair of leaves. At the root node, a root configuration can be obtained by choosing ancestral configurations at each of the cherry nodes and augmenting the configuration with leaves that do not descend from cherry nodes.

Noting $\log(1+x) \leq x$ for $x > 0$, for each ordered unlabeled history $t$ with size $|t| \geq 3$, we have
$$\log \left[1 + \frac{1}{c_r(t)} \right] \leq \frac{1}{c_r(t)} \leq 2^{- \text{ch}(t)}.$$
By taking expectations, we see that Eq.~(\ref{cherr}) implies Eq.~(\ref{teclem}):
$$\mathbb{E}_n\left[ \log \left(1 + \frac{1}{c_r} \right) \right] \leq \mathbb{E}_n[2^{-\text{ch}}].$$

It remains to verify Eq.~(\ref{cherr}). In their Theorem 2, \cite{DisantoAndWiehe13} studied the generating function $F(x,z)$ counting the number of unlabeled histories $t$ of size $n$ with a given number of cherries, where each unlabeled history $t$ is weighted by its probability ${2^{n-1-\text{ch(t)}}}/{(n-1)!}$ under the Yule--Harding distribution:
$$F(x,z) = \sum_{t} \frac{2^{n-1-\text{ch(t)}}}{(n-1)!} x^{\text{ch}(t)} z^n.$$
The sum proceeds over unlabeled histories (``ranked trees'' in \cite{DisantoAndWiehe13}). The coefficient of $x^h z^n$ in $F(x,z)$ gives the probability of $h$ cherries in unlabeled histories of size $n$ under the Yule--Harding distribution, or equivalently, the probability of $h$ cherries in ordered unlabeled histories of size $n$ selected uniformly at random. Hence, the expectation $\mathbb{E}_n[2^{-\text{ch}}]$ is obtained from the coefficient of $z^n$ in $F(\frac{1}{2},z)$. From \cite{DisantoAndWiehe13},
$$F\bigg(\frac{1}{2},z\bigg) = f(z) = \frac{z e^{ z \sqrt{2} } - z}{ (\sqrt{2} - 2) e^{z \sqrt{2}} + 2 + \sqrt{2} }.$$

% \textcolor{red}{
% It remains to verify Eq.~(\ref{cherr}). The expected value of the variable $2^{- \text{ch}}$ over uniform random ordered unlabeled histories of size $n$ can be calculated
% as the $n$th coefficient in the expansion of the generating function
% $$f(z) = \frac{z e^{ z \sqrt{2} } - z}{ (\sqrt{2} - 2) e^{z \sqrt{2}} + 2 + \sqrt{2} },$$
% which is obtained by substituting $x=\frac{1}{2}$ in the function
% $$F(x,z) = \sum_{t} \frac{2^{n-1-\text{ch(t)}}}{(n-1)!} x^{\text{ch}(t)} z^n$$
% appearing in Theorem 2 of \cite{DisantoAndWiehe13}. The sum proceeds over unlabeled histories (``ranked trees'' in \cite{DisantoAndWiehe13}). $F(x,z)$ is the generating
% function counting the number of unlabeled histories $t$ of size $n$ with a given number of cherries, where each unlabeled history $t$ is weighted by its probability
% ${2^{n-1-\text{ch(t)}}}/{(n-1)!}$ under the Yule--Harding distribution. }

By Theorem IV.7 of \cite{FlajoletAndSedgewick09} (see also Section \ref{analyticomb}), $\mathbb{E}_n[2^{-\text{ch}}]$ grows exponentially like $[z^n] f(z) \bowtie \alpha^{-n}$, where $\alpha$ is the dominant singularity of $f(z)$. The value of $\alpha$ is the solution of smallest modulus of the equation $ (\sqrt{2} - 2) e^{z \sqrt{2}} + 2 + \sqrt{2} = 0,$ whose left-hand side is the denominator of $f(z)$. Because
$$\alpha = \frac{1}{\sqrt{2}} \, \log \left(\frac{2+ \sqrt{2}}{2-\sqrt{2}} \right) = \frac{\sqrt{2}\log(3+2\sqrt{2})}{2} \approx 1.246,$$
$\alpha^{-1} \approx 0.802$ and thus, conservatively, $\mathbb{E}_n[2^{-\text{ch}}] = \mathcal{O}(0.9^n)$. Hence, $\mathbb{E}_n[ \log (1 + \frac{1}{c_r} ) ]$ also decays to 0 as $\mathcal{O}(0.9^n)$. $\Box$

\medskip

Considering as in Section \ref{addp} the additive tree parameter $F(t) = \log[c_r(t)+1]$, by Lemma \ref{le} we have demonstrated that the associated toll function $f(t) = \log[1+1/c_r(t)]$ satisfies
\begin{equation}
\label{condth}
\frac{\sum_t f(t)}{F_{n-1}} = \mathbb{E}_n\left[ \log \left(1 + \frac{1}{c_r} \right) \right] = \mathcal{O}(0.9^n),
\end{equation}
where the sum proceeds over all $(n-1)!$ ordered unlabeled histories $t$ of size $n$ (Eq.~(\ref{fact})). Eq.~(\ref{condth}), together with the fact that $f(t)$ is bounded because $c_r(t) \geq 1$ for $|t| \geq 2$, show that the hypotheses of Theorem 4.2 of \cite{Wagner15} are satisfied. By applying the theorem, we can conclude that for an ordered unlabeled history $t$ of size $n$ selected uniformly at random, the standardized version of the random variable $F(t) = \log[c_r(t)+1]$ converges asymptotically to a normal distribution with mean 0 and variance $1$. By the same theorem, the mean and variance of $F(t) = \log[c_r(t)+1]$ grow respectively like $\mu n$ and $\sigma^2 n$, for two constants
\begin{eqnarray}
\label{muu}
\mu &=& \sum_t \frac{2 f(t)}{(|t|+1)!} \approx 0.351 ,
\\\label{sigg} \nonumber
\sigma^2 & = & \sum_t \frac{2 f(t)[2 F(t) - f(t)]}{(|t|+1)!} - \mu^2 + \sum_{t_1} \sum_{t_2} \frac{4 f(t_1) f(t_2)}{(|t_1|+1)!(|t_2|+1)!} \\\nonumber
&& \times \bigg[ \frac{(|t_1|-1)(|t_2|-1)}{|t_1| + |t_2| - 1} - |t_1| - |t_2| + 2 +
\frac{(|t_1|-1)(|t_2|-1)}{(|t_1| + |t_2|)(|t_1| + |t_2| + 1)} \\
&& + \frac{(|t_1|-1)^2(|t_2|-1)^2}{(|t_1| + |t_2| - 1)(|t_1| + |t_2|)(|t_1| + |t_2| + 1)}
\bigg] \approx 0.008.
\end{eqnarray}
Note that the sums in Eqs.~(\ref{muu}) and (\ref{sigg}) are defined over all ordered unlabeled histories, but that the approximations have been calculated by disregarding histories of size strictly larger than 15 and 12 in the sums for $\mu$ and $\sigma^2$, respectively. The equivalence of Lemma \ref{corr} between the distribution of the number of root configurations over uniformly distributed ordered unlabeled histories and the distribution of the number of root configurations over Yule--Harding distributed labeled topologies, coupled with the fact that the difference $\log(c_r+1) - \log c_r = \log(1 + 1/c_r)$ is small, finally yields the following proposition.
\begin{prop}
\label{lognormyul}
The logarithm of the number of root configurations in a labeled topology of size $n$ selected under the Yule--Harding distribution, rescaled as $({\log c_r - \mathbb{E}_n[\log c_r]})/{\sqrt{\mathbb{V}_n[\log c_r]}}$, converges to a standard normal distribution, where $\mathbb{E}_n[\log c_r] \sim \mu n$ and
$\mathbb{V}_n[\log c_r] \sim \sigma^2 n$ for $(\mu, \sigma^2) \approx (0.351, 0.008)$.
\end{prop}

%%%%%%%%%%%%%%%%%%%%%%%%%%%%%%%%%%%%%%%%%%%%%%%%%%%%%%%%%%%%%%
%%%%%%%%%%%%%%%%%%%%%%%% Figure 5 %%%%%%%%%%%%%%%%%%%%%%%%%%%%
%%%%%%%%%%%%%%%%%%%%%%%%%%%%%%%%%%%%%%%%%%%%%%%%%%%%%%%%%%%%%%
\begin{figure}
\begin{center}
\includegraphics*[scale=0.50,trim=0 0 0 0]{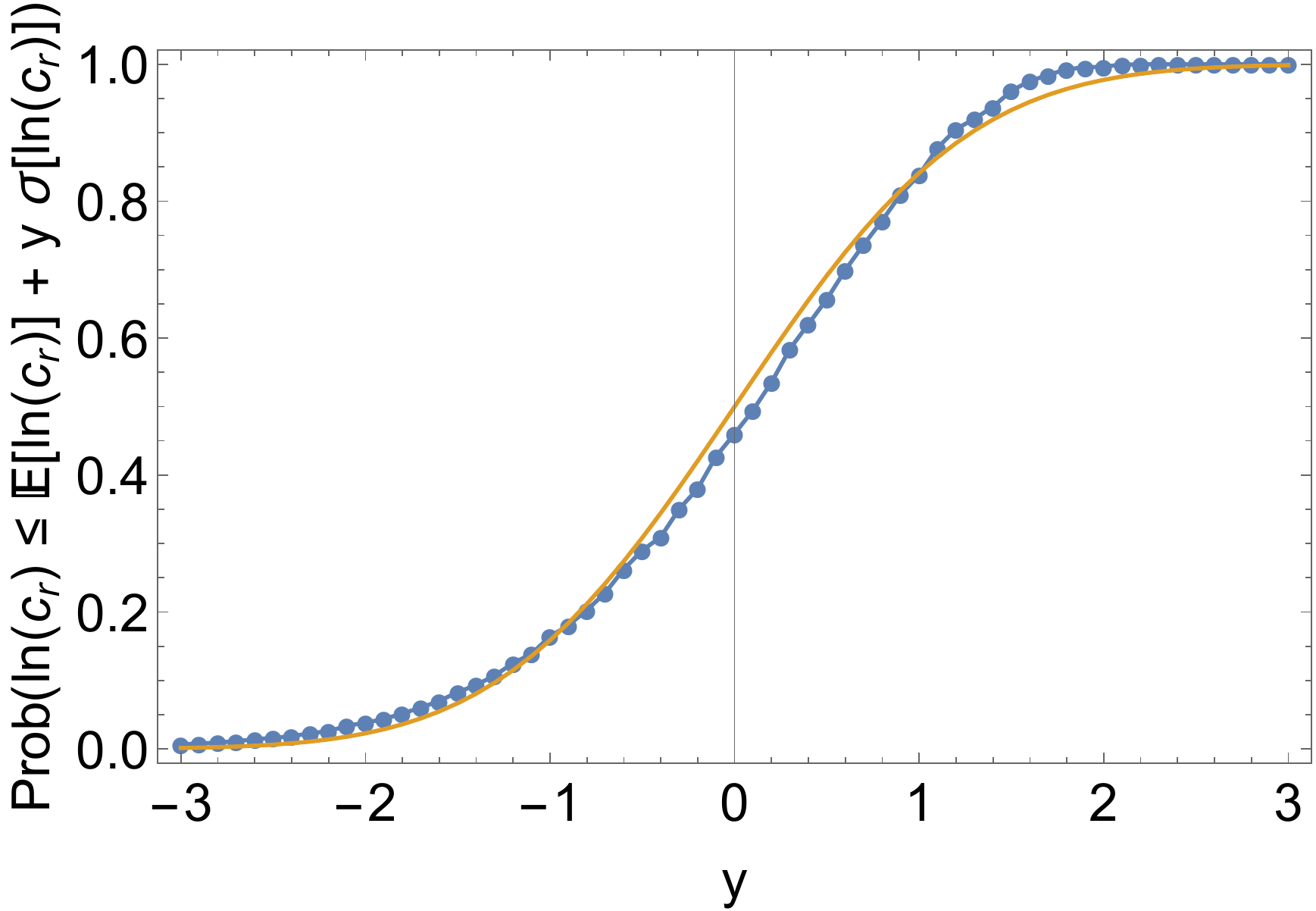}
\end{center}
\vspace{-.7cm}
\caption{{\small Cumulative distribution of the natural logarithm of the number of root configurations for labeled topologies of size $n=15$ considered under the Yule--Harding distribution (dotted line). Each dot has its abscissa determined by a value of $y$ ranging in the interval $y \in [-3,3]$ in steps of $0.1$. Given $y$, the quantity plotted is the probability that a labeled topology with $n=15$ chosen at random under the Yule--Harding distribution has a number of root configurations less than or equal to $\exp\left( \mathbb{E}[\log c_r] + y \sigma[\log c_r] \right)$, where $\mathbb{E}[\log c_r]$ and $\sigma[\log c_r$] are respectively the mean and the standard deviation of the logarithm of the number of root configurations for Yule--Harding distributed labeled topologies of $n=15$ taxa (Proposition \ref{lognormyul}). The solid line is the cumulative distribution of a Gaussian random variable with mean 0 and variance 1.
}}
\label{cumuy}
\end{figure}
%%%%%%%%%%%%%%%%%%%%%%%%%%%%%%%%%%%%%%%%%%%%%%%%%%%%%%%%%%%%%%

For fixed $n$, we can compute the exact distribution of $c_r$ (and $\log c_r$) under the Yule--Harding distribution across all labeled topologies of size $n$ as in Proposition \ref{lemRn}ii. Similarly to the computations in Fig.~\ref{cumuni}, we can weight the counts of root configurations for unlabeled topologies by their Yule--Harding probabilities~\cite[][p.~47]{Steel16}. Fig.~\ref{cumuy} shows the cumulative distribution $\mathbb{P}[ \log c_r \leq \mathbb{E}[\log c_r] + y\sigma[\log c_r] ]$ plotted as a function of $y$, when labeled topologies of size $n = 15$ are selected under the Yule--Harding distribution. The distribution is close to the cumulative distribution of a Gaussian random variable with mean $0$ and variance 1.

%%%%%%%%%%%%%%%%%%%%%%%%%%%%%%%%%%%%%%%%%%%%%%%%%%%%%%%%%%%%%%
%%%%%%%%%%%%%%%%%%%%%%%%%%%%%%%%%%%%%%%%%%%%%%%%%%%%%%%%%%%%%%
%%%%%%%%%%%%%%%%%%%%%%%%%%%%%%%%%%%%%%%%%%%%%%%%%%%%%%%%%%%%%%

\subsection{Mean number of root configurations}
\label{mmy}

In Section \ref{yulelognormal}, we have analyzed distributional properties of the logarithm of the number of root configurations considered over labeled topologies of given size selected under the Yule--Harding distribution. In this section, we study the mean number of root configurations under the Yule--Harding distribution.

From Lemma \ref{corr}, the mean number of root configurations in a random labeled topology of size $n$ selected under the Yule--Harding distribution is also the mean number of root configurations in a uniform random ordered unlabeled history of $n$ taxa. To calculate this mean, we use the distributional recurrence in Proposition \ref{lemRn} for the variable $R_n$ and, by applying generating functions and singularity analysis, we obtain the following result.
\begin{prop}
\label{propken}
The mean number of root configurations in an ordered unlabeled history of size $n$ selected uniformly at random satisfies the asymptotic relation
$\mathbb{E}[R_n] \sim k_e^{n}$, where $k_e = 1/( 1-e^{-2\pi\sqrt{3}/9} )$.
\end{prop}
\noindent \emph{Proof.} Set $e_n \equiv \mathbb{E}[R_n]$. Then $\mathbb{E}[ R_{I_n} \, R^{*}_{n-I_n}] = \sum_{j=1}^{n-1} \mathbb{P}[I_n=j] \, \mathbb{E}[ R_{j} \, R^{*}_{n-j}] = \frac{1}{n-1} \sum_{j=1}^{n-1} \mathbb{E}[ R_{j}] \, \mathbb{E}[R^{*}_{n-j}]$. Proposition \ref{lemRn} yields for $n \geq 2$ the recurrence
\begin{equation}
\label{receqtnen}
e_n=\displaystyle 1 +\dfrac{1}{n-1}\sum_{j=1}^{n-1}e_{j}e_{n-j}+\dfrac{2}{n-1}\sum_{j=1}^{n-1}e_{j},
\end{equation}
with initial condition $e_1=0$.

Defining the generating function
\begin{equation}
\label{expE}
E(z)\equiv\sum_{n=1}^\infty e_n z^n = z^2 + 2z^3 + \frac{10}{3} z^4 + \frac{31}{6} z^5 + \ldots,
\end{equation}
the recurrence in Eq.~(\ref{receqtnen}) translates into the Riccati differential equation
\begin{equation}
\label{equaE1}
z E'(z) = E(z)^2 + \frac{1+z}{1-z} \, E(z) + \frac{z^2}{(1-z)^2},
\end{equation}
with initial condition $E(0) = 0$. To obtain the differential equation, we have multiplied both sides of Eq.~(\ref{receqtnen}) by $(n-1)z^n$, summed for $n \geq 1$, and then used the facts that
$\sum_{n=1}^\infty (n-1) e_n z^n = z E'(z) - E(z)$,
$\sum_{n=1}^\infty (n-1) z^n = z^2 \, [1/(1-z)]' = z^2/(1-z)^2$,
$\sum_{n=1}^\infty ( \sum_{j=1}^{n-1}e_{j}e_{n-j} ) z^n = E(z)^2$, and
$\sum_{n=1}^\infty ( \sum_{j=1}^{n-1} e_{j} ) z^n = E(z) [1/(1-z) - 1]$.

Solving the differential equation yields
\begin{equation}
\label{formE}
E(z) = \frac{2 z \sin\left( \frac{\sqrt{3}}{2} \log(1-z) \right)}{(z-1)\left[ \sqrt{3} \cos\left( \frac{\sqrt{3}}{2} \log(1-z) \right) + \sin\left( \frac{\sqrt{3}}{2} \log(1-z) \right) \right]}.
\end{equation}
In particular, we find that the singularities of $E(z)$ are at $z=1$ and at $z = \alpha \equiv 1-e^{-2\pi\sqrt{3}/9} \approx 0.702$, where the latter is the unique root of the factor
\begin{equation}
\label{factt}
\sqrt{3} \cos \bigg[ \frac{\sqrt{3}}{2} \log(1-z) \bigg] + \sin \bigg[\frac{\sqrt{3}}{2} \log(1-z) \bigg]
\end{equation}
appearing in the denominator of Eq.~(\ref{formE}). The expansion of $E(z)$ at its dominant singularity $z = \alpha$ looks like
\begin{equation*}
E(z) \overset{z \rightarrow \alpha}{\sim} \frac{1}{1- \frac{z}{\alpha}},
\end{equation*}
which can be obtained by plugging the Taylor expansion $-\sqrt{3} e^{+2 \pi \sqrt{3}/9} (z-\alpha)$ of the factor (\ref{factt}) in the denominator of Eq.~(\ref{formE}). By Theorem VI.4 of \cite{FlajoletAndSedgewick09} (see also Section \ref{analyticomb}), we finally obtain
\begin{equation*}
[z^n]E(z) \sim [z^n] \left(\frac{1}{1- \frac{z}{\alpha}} \right) = \alpha^{-n},
\end{equation*}
as $n\rightarrow \infty$. $\Box$

\medskip

The next proposition follows immediately from Proposition \ref{propken}.
\begin{prop}
\label{propyu}
The mean number of root configurations in a labeled topology of size $n$ selected at random under the Yule--Harding distribution has asymptotic growth
$\mathbb{E}_n[c_r] \sim k_e^{n}$, where $k_e = 1/(1-e^{-2\pi\sqrt{3}/9}) \approx 1.42538682$. Furthermore, the mean total number of configurations has asymptotic growth $\mathbb{E}_n[c] \bowtie \mathbb{E}_n[c_r]$.
\end{prop}

%%%%%%%%%%%%%%%%%%%%%%%%%%%%%%%%%%%%%%%%%%%%%%%%%%%%%%%%%%%%%%
%%%%%%%%%%%%%%%%%%%%%%%% Figure 6 %%%%%%%%%%%%%%%%%%%%%%%%%%%%
%%%%%%%%%%%%%%%%%%%%%%%%%%%%%%%%%%%%%%%%%%%%%%%%%%%%%%%%%%%%%%
\begin{figure}
\begin{center}
\includegraphics*[scale=0.50,trim=0 0 0 0]{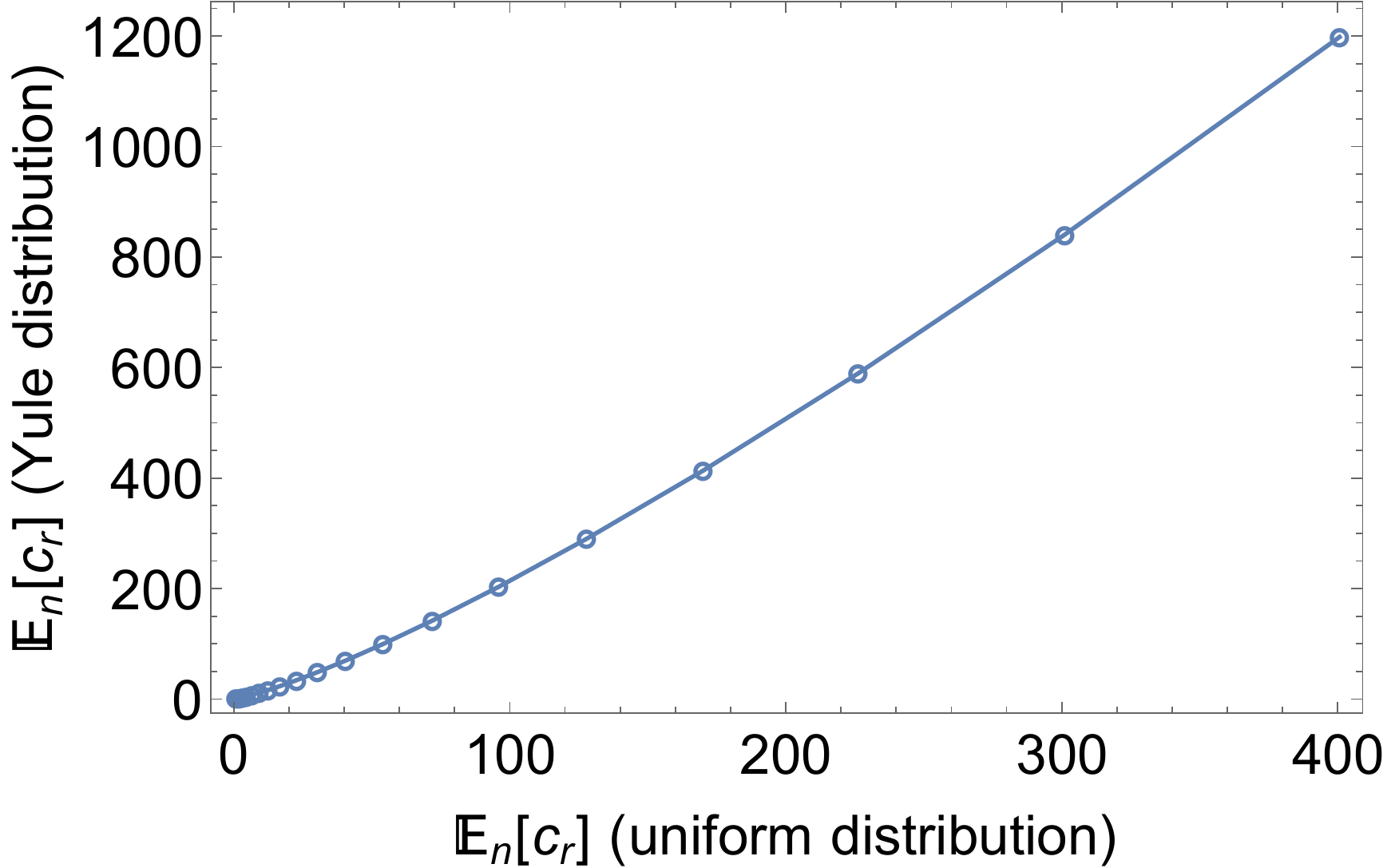}
\end{center}
\vspace{-.7cm}
\caption{{\small Mean number of root configurations of labeled topologies of size $n$ under the Yule--Harding and uniform distributions, for $2 \leq n \leq 20$. Values for the uniform distribution are computed from the power series expansion of Eq.~(33) of \cite{DisantoAndRosenberg17}; values for Yule--Harding are computed from the power series expansion of Eq.~(\ref{formE}).
}}
\label{unifVSyule}
\end{figure}
%%%%%%%%%%%%%%%%%%%%%%%%%%%%%%%%%%%%%%%%%%%%%%%%%%%%%%%%%%%%%%

For small tree size ($n \leq 20$), we plot in Fig.~\ref{unifVSyule} the mean number of root configurations for a random tree of size $n$ selected under the Yule--Harding distribution as a function of the mean number of root configurations under the uniform distribution. The mean is greater for the Yule--Harding distribution, but the two quantities are highly correlated, with Pearson's correlation coefficient approximately 0.995.

%%%%%%%%%%%%%%%%%%%%%%%%%%%%%%%%%%%%%%%%%%%%%%%%%%%%%%%%%%%%%%
%%%%%%%%%%%%%%%%%%%%%%%%%%%%%%%%%%%%%%%%%%%%%%%%%%%%%%%%%%%%%%
%%%%%%%%%%%%%%%%%%%%%%%%%%%%%%%%%%%%%%%%%%%%%%%%%%%%%%%%%%%%%%

\subsection{Variance of the number of root configurations}
\label{vvy}

In this section, we analyze the asymptotic growth of the variance of the number of root configurations under the Yule--Harding distribution. In particular, by using Lemma \ref{corr}, we study the variance of the number of root configurations in a uniform random ordered unlabeled history of size $n$.

Following Section \ref{mmy} and squaring Eq.~(\ref{eqRn}), we obtain a recurrence for $s_n \equiv \mathbb{E}[R_n^2]$.
For $n \geq 2$,
\begin{equation}
\label{recurss}
s_n = 1 + \frac{1}{n-1} \sum_{j=1}^{n-1} s_j \, s_{n-j}
+ \frac{2}{n-1} \sum_{j=1}^{n-1} s_j
+ \frac{4}{n-1} \sum_{j=1}^{n-1} s_j \, e_{n-j}
+ \frac{4}{n-1} \sum_{j=1}^{n-1} e_j \, e_{n-j}
+ \frac{4}{n-1} \sum_{j=1}^{n-1} e_j ,
\end{equation}
with initial condition $s_1 = 0$.

Starting from this recurrence, a symbolic calculation similar to that used to derive Eq.~(\ref{equaE1}) shows that the generating function $S(z) \equiv \sum_{n=1}^\infty s_n z^n = z^2 + 4z^3 + \frac{34}{3}z^4 + \frac{55}{2}z^5 \ldots$
satisfies the Riccati differential equation
\begin{equation}
z \, S'(z) = S(z)^2 - S(z) \left[ \frac{1+z}{z-1} - 4 E(z) \right] + \frac{[ z - 2 (z-1) E(z) ]^2}{(z-1)^2}.
\end{equation}
This equation can be written
\begin{equation}
\label{eqs}
S'(z) = g_2(z) \, S(z)^2 + g_1(z) \, S(z) + g_0(z)
\end{equation}
by setting
$$\bigg(g_2(z),g_1(z),g_0(z)\bigg) \equiv \left( \frac{1}{z} , \left( 4 E(z) - \frac{1+z}{z-1} \right) \frac{1}{z} , \frac{[ z - 2 (z-1) E(z)]^2}{z (z-1)^2} \right).$$

By substituting $U(z) \equiv \exp[\int_0^z S(x)/(-x) \, dx]$, we obtain $S(z) =
%- \frac{U'(z)}{g_2 U(z)} =
{-z U'(z)}/{U(z)},$ and Eq.~(\ref{eqs}) can be rewritten as a second-order linear differential equation equation
\begin{equation}
\label{diffeqtn}
U''(z) - \left(g_1(z) + \frac{g_2'(z)}{g_2(z)} \right) U'(z) + g_2(z) \, g_0(z) \, U(z) = 0.
\end{equation}
%with conditions $U(0)=1$ and $U'(0)=0$ derived from $U = exp[\int -g_2 \, S]$.
The coefficients of Eq.~(\ref{diffeqtn}) are analytic functions for $|z| < 0.702$, with a removable singularity at $z=0$ as the expansion (\ref{expE}) of $E(z)$ starts with a quadratic non-zero term. Using existence results for the solutions of second-order ordinary differential equations, $U(z)$ must be analytic for $|z| < 0.702$, the constant being the radius of convergence of $E(z)$ as determined in the proof of Proposition \ref{propken}. Therefore, also $U'(z)$ is analytic for $|z| < 0.702$, and thus $S(z)$ is a meromorphic function on this domain, being a quotient of two analytic functions. To analyze the singularities of a meromorphic function, one must locate the possible roots of its denominator function. In our case, the set of singularities of $S(z)$ consists of the roots of $U(z)$. In particular, by studying in the Appendix the function $U(z)$ in $\mathcal{B}\equiv \{z\in\mathbb{C}\ :\ |z|\leq \frac{1}{2} \}$, we find that $S(z)$ has a unique dominant singularity $\alpha \approx 0.4889986317$, the unique and simple root of $U(z)$ within $\mathcal{B}$ (Proposition \ref{rootOfU}).

As a consequence, we can write $U(z) = (z - \alpha) \tilde{U}(z)$, with $\tilde{U}(\alpha) \neq 0$ and $U'(\alpha) = (- \alpha) \tilde{U}(\alpha) \neq 0$. Therefore, for $ z \rightarrow \alpha$ the generating function $S(z)$ admits the expansion
\begin{equation*}
%\label{asis}
S(z) = \frac{-z U'(z)}{U(z)} \overset{z \rightarrow \alpha}{\sim} \frac{(-\alpha) [U'(\alpha) + U''(\alpha)(z-\alpha) + \ldots]}{U(\alpha) + U'(\alpha)(z-\alpha) + \ldots} \overset{z \rightarrow \alpha}{\sim} \frac{(-\alpha)U'(\alpha)}{U'(\alpha)(z-\alpha)} = \frac{- \alpha}{z-\alpha} = \frac{1}{1 - \frac{z}{\alpha}}.
\end{equation*}
From Theorem VI.4 of \cite{FlajoletAndSedgewick09} (see also Section \ref{analyticomb}), we can thus recover the asymptotic growth of the associated coefficients
\begin{equation}
\label{asyR2}
\mathbb{E}[R_n^2] = [z^n] S(z) \sim [z^n] \left(\frac{1}{1 - \frac{z}{\alpha}} \right) = \alpha^{-n},
\end{equation}
and hence derive the asymptotic growth of the variance $\mathbb{V}[R_n]$. In particular, we have the following result.
\begin{prop}
\label{propyuv}
The variance of the number of root configurations in a labeled topology of size $n$ selected at random under the Yule--Harding distribution has asymptotic growth $\mathbb{V}_n[c_r] \sim k_v^n$, where $k_v \approx 2.0449954971$. Furthermore, the variance of the total number of configurations has asymptotic growth $\mathbb{V}_n[c] \bowtie \mathbb{V}_n[c_r]$.
\end{prop}

\noindent \emph{Proof.} For uniformly distributed ordered unlabeled histories of size $n$, Eq.~(\ref{asyR2}) yields $\mathbb{E}[R_n^2] \sim k_v^{n}$, $k_v \equiv 1/\alpha \approx 2.0449954971$. From Proposition \ref{propken}, $\mathbb{E}[R_n]^2 \sim (k_e^2)^{n}$, with $k_e^2 \approx 2.03$. Because $k_v>k_e^2$, as $n\rightarrow\infty$ we obtain
$${\mathbb V}[R_n] = \mathbb{E}[R_n^2] - \mathbb{E}[R_n]^2 \sim k_v^n.$$
By Lemma \ref{corr}, the variance of the variable $R_n$ is the variance of the number of root configurations considered over labeled topologies of $n$ taxa selected under the Yule--Harding distribution. $\Box$

\medskip

%%%%%%%%%%%%%%%%%%%%%%%%%%%%%%%%%%%%%%%%%%%%%%%%%%%%%%%%%%%%%%
%%%%%%%%%%%%%%%%%%%%%%%% Figure 7 %%%%%%%%%%%%%%%%%%%%%%%%%%%%
%%%%%%%%%%%%%%%%%%%%%%%%%%%%%%%%%%%%%%%%%%%%%%%%%%%%%%%%%%%%%%
\begin{figure}
\begin{center}
\includegraphics*[scale=0.50,trim=0 0 0 0]{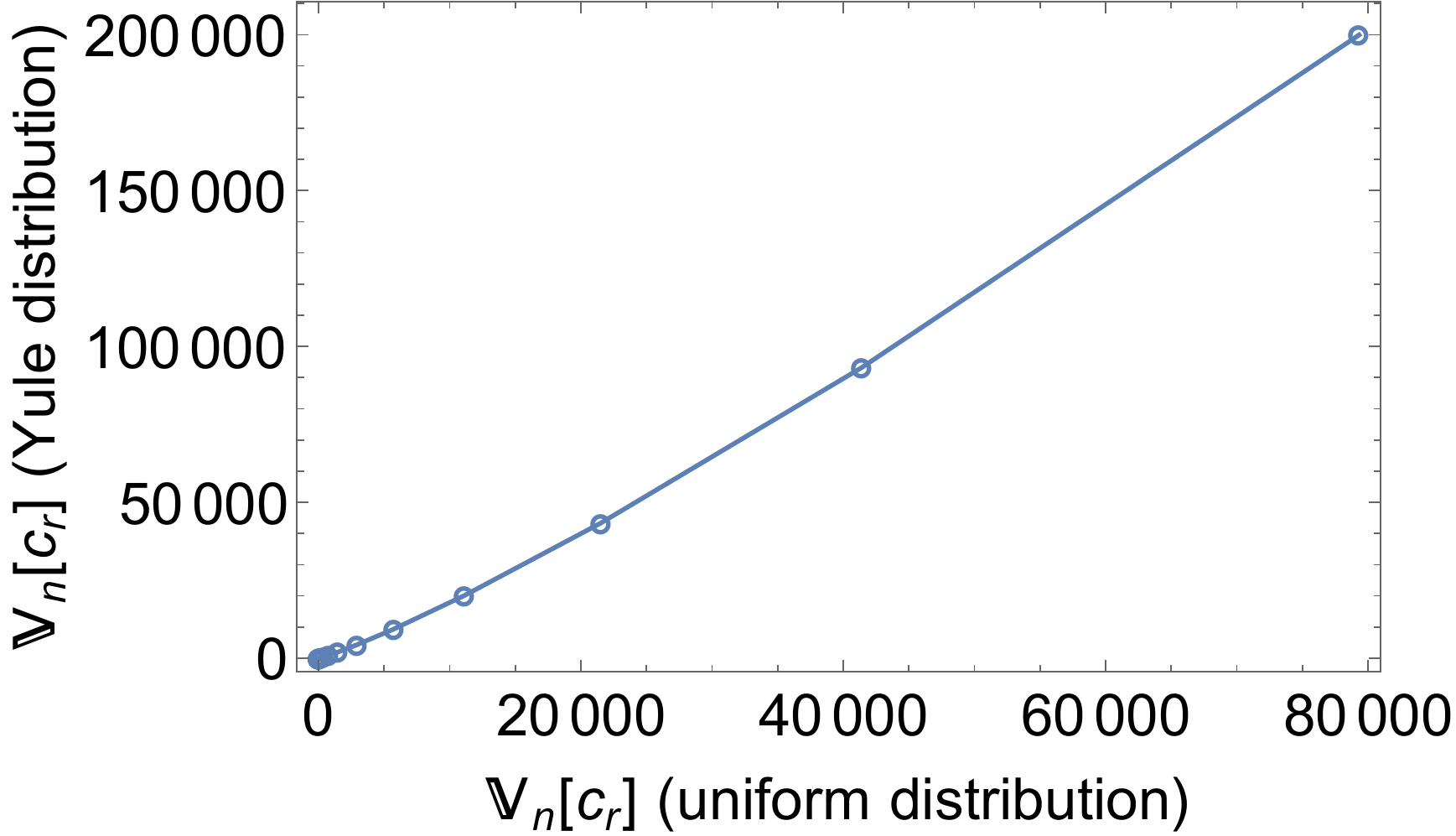}
\end{center}
\vspace{-.7cm}
\caption{{\small Variance of the number of root configurations of labeled topologies of size $n$ under the Yule--Harding and uniform distributions, for $2 \leq n \leq 20$. Values for the uniform distribution are computed from the power series expansion of Eq.~(39) of \cite{DisantoAndRosenberg17}; values for Yule--Harding are computed from Eqs.~(\ref{recurss}) and (\ref{receqtnen}).
}}
\label{unifVSyulevar}
\end{figure}
%%%%%%%%%%%%%%%%%%%%%%%%%%%%%%%%%%%%%%%%%%%%%%%%%%%%%%%%%%%%%%

For small tree size ($n \leq 20$), we plot in Fig.~\ref{unifVSyulevar} the variance of the number of root configurations for a random tree of size $n$ selected under the Yule--Harding distribution as a function of the variance of the number of root configurations for a random uniform tree of the same size. As was true of the mean, the Yule--Harding and uniform distributions on labeled topologies give correlated variances (correlation coefficient 0.997).

\section{Discussion}

Considering gene trees and species trees with a matching labeled topology $G=S=t$, we have studied distributional properties of the number $c_r$ of root ancestral configurations for labeled topologies $t$ of fixed size under two probability models, the uniform model and the Yule--Harding model (Table \ref{risfin}). We have made use of techniques of analytic combinatorics, relying on equivalences across tree types (Section \ref{secEquivalences}), and making particular use of results of \cite{Wagner15} on distributional properties of additive tree parameters for several families of trees.

Extending results of \cite{DisantoAndRosenberg17}, for the uniform model we have shown that the logarithm of the number of root configurations, when standardized, converges asymptotically to a standard normal distribution (Proposition \ref{lognormunil}). Under the Yule--Harding distribution, as is the case for uniformly distributed labeled topologies, the logarithm of the number of root configurations, when standardized, converges to a standard normal distribution (Proposition \ref{lognormyul}). We have also determined the asymptotic growth of the mean and the variance of the number of root configurations, finding that under the Yule--Harding model, $\mathbb{E}_n[c_r] \sim 1.425^n$ (Proposition \ref{propyu}) and $\mathbb{V}_n[c_r] \sim 2.045^n$ (Proposition \ref{propyuv}). As $\mathbb{E}_n[c] \bowtie \mathbb{E}_n[c_r]$ and $\mathbb{V}_n[c] \bowtie \mathbb{V}_n[c_r]$, we also recover the exponential growth rate of the mean and the variance of the total number of configurations under the Yule--Harding model.

\begin{figure}
\begin{center}
\includegraphics*[scale=0.36,trim=0 0 0 0]{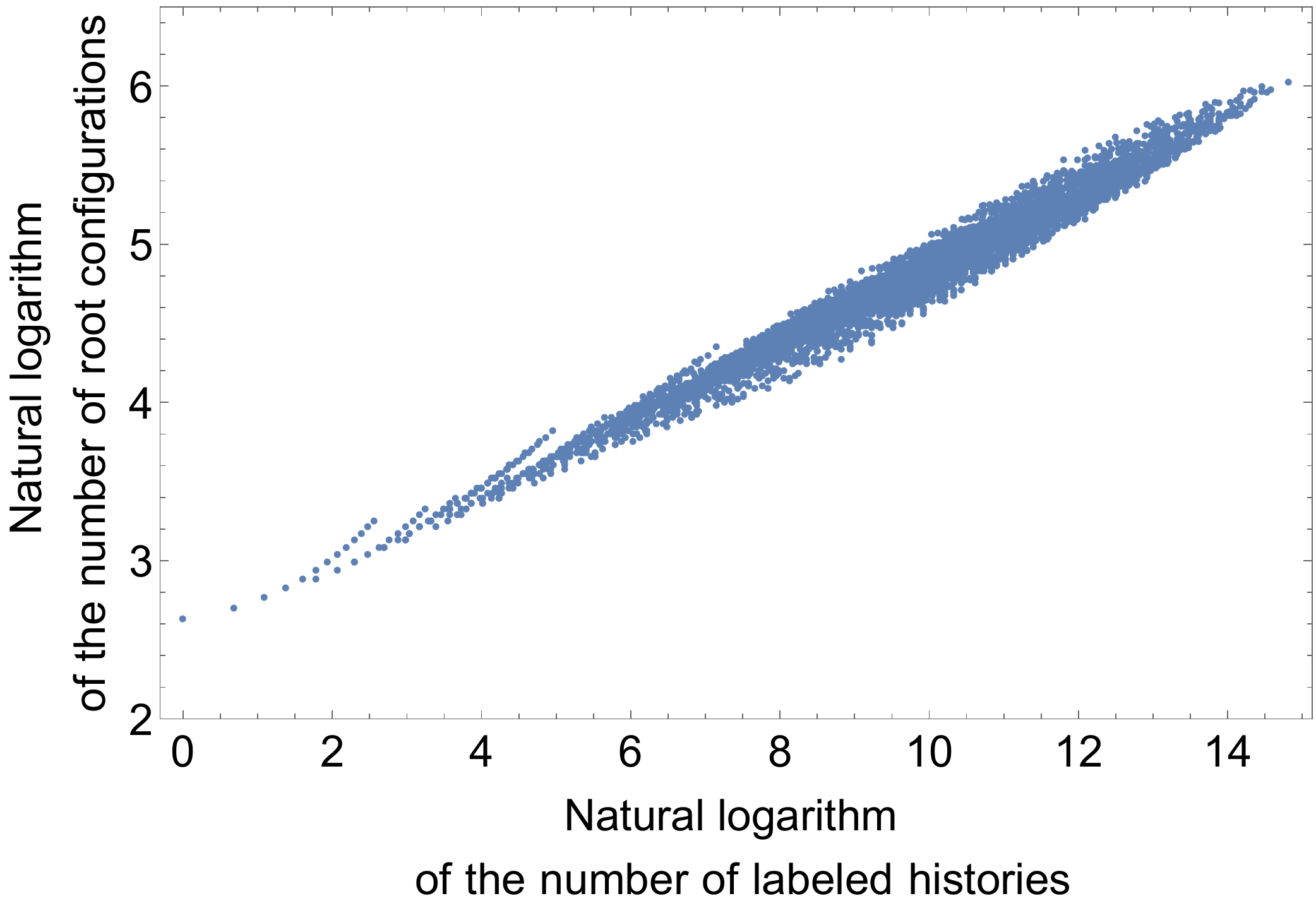}
\end{center}
\vspace{-.7cm}
\caption{{\small Natural logarithm of the number of root configurations and natural logarithm of the number of labeled histories for a representative labeling of each unlabeled topology of size $n=15$. The number of points plotted is 4850, the number of unlabeled topologies with $n=15$ taxa. The Pearson correlation is approximately 0.987 (0.784 without log scaling).
}}
\label{figbal}
\end{figure}

%%%%%%%%%%%%%%%%%%%%%%%%%%%%%%%%%%%%%%%%%%%%%%%%%%%%%%%%%%%%%%
%%%%%%%%%%%%%%%%%%%%%%%% Table 1 %%%%%%%%%%%%%%%%%%%%%%%%%%%%%
%%%%%%%%%%%%%%%%%%%%%%%%%%%%%%%%%%%%%%%%%%%%%%%%%%%%%%%%%%%%%%
\begin{table}
\caption{Distributional properties of the number of root and total configurations.}
\begin{center}
\label{risfin}
\vskip -.7cm
\fontsize{9}{11}\selectfont
\begin{tabular}{ ||p{2.4cm}|p{2.cm}||c c|c c|| }
\hline \hline
\multicolumn{2} {||c||}{Results}& \multicolumn{2} {c|}{Uniform model} &\multicolumn{2} {c|}{Yule--Harding model} \\ \hline \hline
\multirow{4}{*}{ \parbox{2.4cm}{Root \\ configurations} }
& Mean     & $\mathbb{E}_n[c_r] \sim 1.225 \cdot 1.333^n$ & Eq.~(\ref{meandis}) & $\mathbb{E}_n[c_r] \sim 1.425^n$ & Proposition \ref{propyu}  \\
& Variance & $\mathbb{V}_n[c_r] \sim 1.405 \cdot 1.822^n$ & Eq.~(\ref{vardis})  & $\mathbb{V}_n[c_r] \sim 2.045^n$ & Proposition \ref{propyuv} \\
& \multirow{2}{*}{\parbox{2.cm}{Lognormal \\ distribution}} & $\mathbb{E}_n[\log c_r] \sim 0.272 \cdot n$ & Proposition \ref{lognormunil} & $\mathbb{E}_n[\log c_r ] \sim 0.351 \cdot n$ & Proposition \ref{lognormyul} \\
&                                                           & $\mathbb{V}_n[\log c_r] \sim 0.034 \cdot n$ & Proposition \ref{lognormunil} & $\mathbb{V}_n[\log c_r] \sim 0.008 \cdot n$ & Proposition \ref{lognormyul} \\ \hline
\multirow{2}{*}{ \parbox{2.4cm}{Total \\ configurations}}
& Mean     & $\mathbb{E}_n[c] \bowtie 1.333^n$ & Eq.~(\ref{meandiss}) & $\mathbb{E}_n[c] \bowtie 1.425^n$ & Proposition \ref{propyu} \\
& Variance & $\mathbb{V}_n[c] \bowtie 1.822^n$ & Eq.~(\ref{vardiss})  & $\mathbb{V}_n[c] \bowtie 2.045^n$ & Proposition \ref{propyuv} \\
\hline
\hline
\end{tabular}
\end{center}
\end{table}
%%%%%%%%%%%%%%%%%%%%%%%%%%%%%%%%%%%%%%%%%%%%%%%%%%%%%%%%%%%%%%
The difference in results for the uniform and Yule--Harding models, along with the results of \cite{DisantoAndRosenberg17}, suggests a role for tree balance in predicting the number of root configurations. By considering a representative labeling for each unlabeled topology of size $n=15$, in Figure \ref{figbal} we plot on a logarithmic scale the number of root configurations as a function of the number of labeled histories, the latter calculated as in Eq.~(\ref{numblab}). The figure shows that the two quantities are correlated: highly balanced labeled topologies---which tend to have a larger number of labeled histories (Section \ref{yhdistr})---in general have a larger number of root configurations.

\begin{figure}
\begin{center}
\includegraphics*[scale=0.86,trim=0 0 0 0]{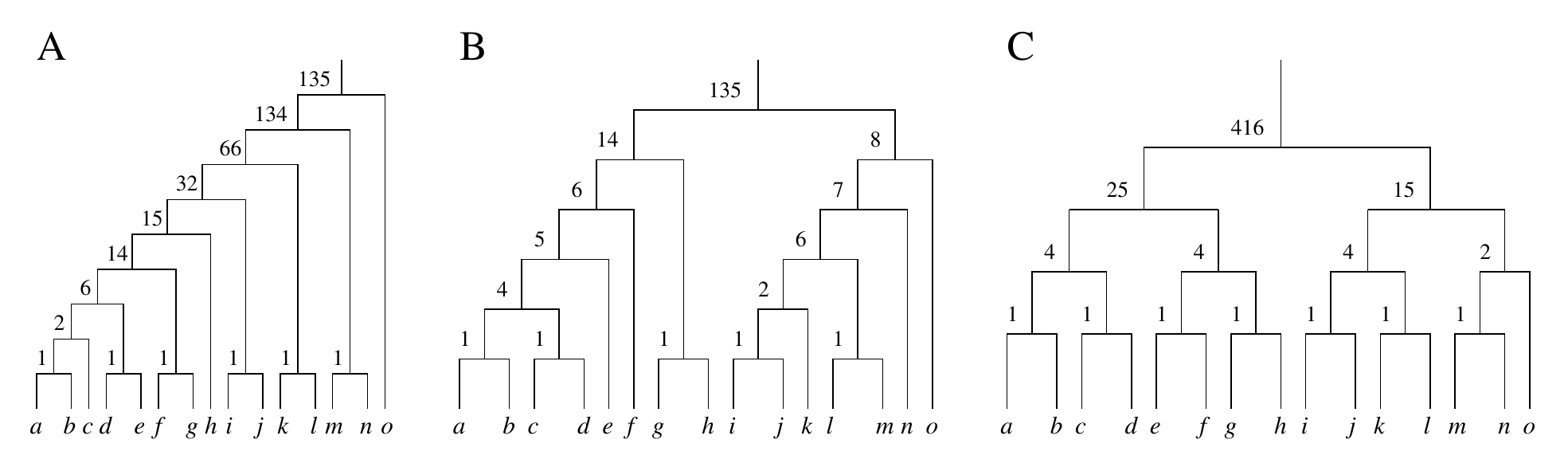}
\end{center}
\vspace{-.7cm}
\caption{{\small The number of ancestral configurations at the internal nodes of three labeled topologies of size $n=15$. {\bf (A, B)} Two labeled topologies in which the number of root configurations is the mean number $c_r = 135$ of root configurations calculated across the set of representative labelings of the unlabeled topologies of size $15$. In this set, the labeled topologies in (A) and (B) have respectively the largest number 61776 and smallest number 14400 of labeled histories. {\bf (C)} The labeled topology with 15 taxa that has the most root configurations (416) and the most labeled histories (2745600).
}}
\label{figtreesl}
\end{figure}

In particular, the largest number of root configurations is possessed by the balanced labeled topology depicted in Figure \ref{figtreesl}C, which also has the largest number of labeled histories, 2745600. The trend in this example is confirmed by our asymptotic results. Under the Yule--Harding probability model, which gives more weight to balanced labeled topologies than does the uniform model, the mean number of root configurations and the mean total number of configurations grow exponentially faster than under the uniform distribution (Table \ref{risfin}). This differing behavior also accords with the proof of \cite{DisantoAndRosenberg17} that balanced and caterpillar trees respectively possess the largest and smallest numbers of root configurations for fixed tree size (Section \ref{ignores}).

%%%%%%%%%%%%%%%%%%%%%%%%%%%%%%%%%%%%%%%%%%%%%%%%%%%%%%%%%%%%%%
%%%%%%%%%%%%%%%%%%%%%%%% Figure 8 %%%%%%%%%%%%%%%%%%%%%%%%%%%%
%%%%%%%%%%%%%%%%%%%%%%%%%%%%%%%%%%%%%%%%%%%%%%%%%%%%%%%%%%%%%%
%%%%%%%%%%%%%%%%%%%%%%%%%%%%%%%%%%%%%%%%%%%%%%%%%%%%%%%%%%%%%%

%%%%%%%%%%%%%%%%%%%%%%%%%%%%%%%%%%%%%%%%%%%%%%%%%%%%%%%%%%%%%%
%%%%%%%%%%%%%%%%%%%%%%%% Figure 9 %%%%%%%%%%%%%%%%%%%%%%%%%%%%
%%%%%%%%%%%%%%%%%%%%%%%%%%%%%%%%%%%%%%%%%%%%%%%%%%%%%%%%%%%%%%
%%%%%%%%%%%%%%%%%%%%%%%%%%%%%%%%%%%%%%%%%%%%%%%%%%%%%%%%%%%%%%%
Several directions naturally arise from our work. First, we focused on root rather than total configurations; although some results for total configurations follow quickly (Table \ref{risfin}), we did not consider total configurations in detail. Second, we assumed that the gene tree and species tree had the same labeled topology, and we did not study nonmatching gene trees and species trees. The nonmatching case merits further analysis, as a nonmatching gene tree labeled topology can have more root and total configurations than the topology that matches the species tree \citep{DisantoAndRosenberg17}. Third, ancestral configurations can be considered up to an equivalence relationship that accounts for symmetries in gene trees \citep{Wu12}. The resulting equivalence classes---the nonequivalent ancestral configurations---are used for calculating probabilities of gene trees in STELLS \citep{Wu12}, with computational complexity that depends on the number of these classes. Some investigation of this number has been carried out by \cite{DisantoAndRosenberg19:bmb} for uniformly distributed matching gene trees and species trees. It would be of interest to see whether the techniques we have used could derive distributional properties of the number of nonequivalent ancestral configurations under the uniform and Yule--Harding probability models.

%%%%%%%%%%%%%%%%%%%%%%%%%%%%%%%%%%%%%%%%%%%%%%%%%%%%%%%%%%%%%%
%%%%%%%%%%%%%%%%%%%%%%%%%%%%%%%%%%%%%%%%%%%%%%%%%%%%%%%%%%%%%%
%%%%%%%%%%%%%%%%%%%%%%%%%%%%%%%%%%%%%%%%%%%%%%%%%%%%%%%%%%%%%%

\subsection*{Appendix. The function $U(z)$ has a unique and simple root of smallest modulus}

In this appendix, we prove that the function $U(z) \equiv \sum_{n=0}^\infty u_n z^n$, which is analytic in the region $|z| < 0.702 $ and there satisfies the differential equation in Eq.~(\ref{diffeqtn}), has a unique and simple root $\alpha$ of smallest modulus. We also calculate the first ten digits of $\alpha \approx 0.4889986317$.

We start in Lemma \ref{lem6} by providing a recurrence for $u_n$, which is then used to find an upper bound of $|u_n|$ in Lemma \ref{uuu}. Next, we consider the set $\mathcal{B} \equiv\{z\in\mathbb{C}\ :\ |z|\leq \frac{1}{2} \}$ in the complex plane and decompose $U(z)$ into a sum $U(z) = U_1(z) + U_2(z)$, where $U_1(z) = \sum_{n = 0}^{100} u_n z^n$ is a polynomial and $U_2(z) =\sum_{n=101}^\infty u_n z^n $. The bound for $|u_n|$ in Lemma \ref{uuu} yields a bound for $|U_1(z)|$ (Lemma \ref{lower-bound}), which in turn implies that $|U_1(z)| > |U_2(z)|$ if $z \in \partial \mathcal{B}$. Hence, by Rouch\'e's theorem we have that inside $\mathcal{B}$, the function $U(z)$ has the same number of roots---considered with their multiplicity---as the polynomial $U_1(z)$. Lemma \ref{lem10} shows that $U_1(z)$ has a unique and simple root inside $\mathcal{B}$, and in Proposition \ref{rootOfU} we conclude the proof of our claim by finding an approximation of the unique and simple root $\alpha$ of $U(z)$ inside $\mathcal{B}$---which turns out to be very close to the root of $U_1(z)$ inside $\mathcal{B}$.

In $U(z) = \sum_{n=0}^\infty u_n z^n$, we have $u_n \equiv [z^n] U(z)$. From Eq.~(\ref{diffeqtn}), we derive a recurrence for $u_n$. Recall that $e_n$ gives the mean number of root configurations in an ordered unlabeled history of size $n \geq 1$.
\begin{lemm}
\label{lem6}
For $n \geq 2$, we have
\begin{equation}
u_n=\dfrac{1}{n(n-1)}\sum_{k=0}^{n-1}(3n-k-3)u_k-\dfrac{4}{n(n-1)}\sum_{k=0}^{n-1}(n-2k-1)e_{n-k}u_k+\dfrac{4}{n(n-1)}\sum_{k=0}^{n-1}\bigg(\sum_{j=0}^{n-k-1}e_j\bigg)u_k,
\label{receqtnun}
\end{equation}
with $u_0=1$ and $u_1=0$.
\end{lemm}
\noindent \emph{Proof.} First notice that for $n \geq 0$, the coefficient of $z^n$ in each term of Eq.~(\ref{diffeqtn}) can be written as
\begin{equation*}
[z^n]U''(z) = (n+2)(n+1)u_{n+2}
\end{equation*}
\begin{equation*}
-[z^n]\left(g_1+\dfrac{g_2'}{g_2}\right)U'(z) = -\sum_{k=0}^{n}(n-k+1)(4e_{k+1}+2)u_{n-k+1}
\end{equation*}
\begin{equation*}
[z^n] g_2 g_0 U(z) = \sum_{k=0}^{n}\bigg[(k+1)+4\sum_{j=0}^{k}e_{j+1}+4\sum_{j=0}^{k+2}e_je_{k-j+2}\bigg]u_{n-k},
\end{equation*}
where for convenience we set $e_0=0$.

Making a substitution to the index of summation, we have
\begin{align*}
&-4\sum_{k=0}^n(n-k+1)e_{k+1}u_{n-k+1}=-4\sum_{k=0}^{n+1}ke_{n-k+2}u_k.
\end{align*}
Hence, the sum for $- [z^n](g_1+{g_2'}/{g_2})U'(z)$ can be simplified as
\begin{align*}
&-[z^n]\left(g_1+\dfrac{g_2'}{g_2}\right)U'(z)=- 4\sum_{k=0}^{n+1}ke_{n-k+2}u_k - 2\sum_{k=0}^n(n-k+1)u_{n-k+1}.
\end{align*}
The second sum in this equation together with the first sum $\sum_{k=0}^n (k+1) u_{n-k}$ of $[z^n]g_2g_0U(z)$ give
\begin{align*}
&-2\sum_{k=0}^n(n-k+1)u_{n-k+1}+\sum_{k=0}^n(k+1)u_{n-k}=\sum_{k=0}^{n+1}(n-3k+1)u_k.
\end{align*}
Furthermore, by setting $n = k+2$ in Eq.~(\ref{receqtnen}), the inner sums of $[z^n]g_2g_0U(z)$ can be rewritten as
\begin{align*}
&4\sum_{j=0}^{k}e_{j+1}+4\sum_{j=0}^{k+1}e_je_{k-j+2}=4(k+1)e_{k+2}-4(k+1)-4\sum_{j=1}^{k+1}e_j.
\end{align*}
Hence, the coefficient of $z^n$ in Eq.~(\ref{diffeqtn}) becomes
$$(n+2)(n+1)u_{n+2} - 4\sum_{k=0}^{n+1}ke_{n-k+2}u_k + \sum_{k=0}^{n+1}(n-3k+1)u_k + \sum_{k=0}^n \bigg[ 4(k+1)e_{k+2}-4(k+1)-4\sum_{j=1}^{k+1}e_j \bigg] u_{n-k}.$$
In this expression, we make two substitutions:
\begin{eqnarray}
\sum_{k=0}^{n}4(k+1)e_{k+2}u_{n-k} & = & \sum_{k=0}^{n+1}4(n-k+1)e_{n-k+2}u_k \\
\sum_{k=0}^{n+1}(n-3k+1)u_k-4\sum_{k=0}^{n}(k+1)u_{n-k}
& = & \sum_{k=0}^{n+1}(n-3k+1)u_k-4\sum_{k=0}^{n}(n-k+1)u_{k} =\sum_{k=0}^{n+1}(-3n+k-3)u_k,
\end{eqnarray}
obtaining
$$(n+2)(n+1)u_{n+2} - 4\sum_{k=0}^{n+1}ke_{n-k+2}u_k + \sum_{k=0}^{n+1}4(n-k+1)e_{n-k+2}u_k + \sum_{k=0}^{n+1}(-3n+k-3)u_k + \sum_{k=0}^n \bigg( -4\sum_{j=1}^{k+1}e_j \bigg) u_{n-k},$$
and thus
$$(n+2)(n+1)u_{n+2} + \sum_{k=0}^{n+1}4(n-2k+1)e_{n-k+2}u_k + \sum_{k=0}^{n+1}(-3n+k-3)u_k + \sum_{k=0}^n \bigg( -4\sum_{j=1}^{k+1}e_j \bigg) u_{n-k}.$$

Finally, because $e_0 = 0$, in this expression we can substitute
\begin{align*}
& \sum_{k=0}^n \bigg( -4\sum_{j=1}^{k+1}e_j \bigg) u_{n-k} = \sum_{k=0}^n \bigg( -4\sum_{j=0}^{k+1}e_j \bigg) u_{n-k} = \sum_{k=0}^n \bigg( -4\sum_{j=0}^{n-k+1}e_j \bigg) u_{k} = \sum_{k=0}^{n+1} \bigg( -4\sum_{j=0}^{n-k+1}e_j \bigg) u_{k},
\end{align*}
obtaining for $n \geq 0$
$$(n+2)(n+1)u_{n+2} + \sum_{k=0}^{n+1}4(n-2k+1)e_{n-k+2}u_k -\sum_{k=0}^{n+1}(3n-k+3)u_k - 4\sum_{k=0}^{n+1}\bigg(\sum_{j=0}^{n-k+1}e_j\bigg)u_k=0,$$
which rescaled is recurrence (\ref{receqtnun}). The starting conditions $u_0 = 1$ and $u_1 = 0$, follow from the fact that $U(0)=1$ and $U'(0)=0$ as $U(z) = \exp[\int_0^z S(x)/(-x) \, dx]$. $\Box$

\medskip

In Lemma \ref{uuu}, we use the recurrence to find an upper bound for $|u_n|$. First, we need an upper bound for $e_n$.
\begin{lemm}
\label{lemEee}
For $n\geq 0$, we have $e_n\leq (\frac{9}{10})(\frac{3}{2})^n$.
\end{lemm}

\noindent \emph{Proof.}
Using the recurrence (\ref{receqtnen}), with the help of computing software we have shown that the inequality holds for $0\leq n\leq 41$. We proceed by induction. Suppose the inequality holds for all $k<n$ with $n>41$. By Eq.~(\ref{receqtnen}),
\begin{align*}
e_n&\leq\displaystyle 1+\dfrac{81}{100(n-1)}\sum_{j=1}^{n-1}\left(\dfrac{3}{2}\right)^n+\dfrac{9}{5(n-1)}\sum_{j=1}^{n-1}\left(\dfrac{3}{2}\right)^j\\
&=1+\dfrac{81}{100}\left(\dfrac{3}{2}\right)^n+\dfrac{18}{5(n-1)}\left(\dfrac{3}{2}\right)^n-\dfrac{27}{5(n-1)}\\
&=\dfrac{9}{10}\left(\dfrac{3}{2}\right)^n-\dfrac{9}{10}\left(\dfrac{1}{10}-\dfrac{4}{n-1}\right)\left(\dfrac{3}{2}\right)^n-\dfrac{27}{5(n-1)}+1.
\end{align*}
In the last step, we can see that a positive number is subtracted from $\frac{9}{10}(\frac{3}{2})^n$ for $n>41$, as $$\dfrac{9}{10}\left(\frac{1}{10}-\dfrac{4}{n-1}\right)\left(\dfrac{3}{2}\right)^n+\dfrac{27}{5(n-1)} - 1 > \frac{9}{10} \frac{1}{400} \left(\dfrac{3}{2}\right)^{42} - 1 > 0.$$
Thus, the claim is proved. $\Box$

\begin{lemm}
\label{uuu}
For $n\geq 0$, we have $|u_n|\leq (\frac{9}{5})^{n}$.
\end{lemm}

\noindent \emph{Proof.}
Using recurrence (\ref{receqtnun}), computing software verifies the inequality for $0\leq n\leq 25$. We proceed by induction. Suppose that the inequality holds for all $k<n$ with $n>25$. For simplicity of computation, instead of the bound in Lemma \ref{lemEee}, we use the more conservative $(\frac{3}{2})^n$ as a bound for $e_n$. With Eq.~(\ref{receqtnun}), we get
\begin{align*}
|u_n|&\leq\displaystyle\dfrac{3}{n}\sum_{k=0}^{n-1}\left(\dfrac{9}{5}\right)^k+\dfrac{4}{n}\sum_{k=0}^{n-1}\left(\dfrac{3}{2}\right)^{n-k}\left(\dfrac{9}{5}\right)^k
+\dfrac{4}{n(n-1)}\sum_{k=0}^{n-1}\left(\sum_{j=0}^{n-k-1}\left(\dfrac{3}{2}\right)^j\right)\left(\dfrac{9}{5}\right)^k\\
&=\displaystyle\dfrac{15}{4n}\left(\dfrac{9}{5}\right)^n-\dfrac{15}{4n}+\dfrac{20}{n}\left(\dfrac{9}{5}\right)^n-\dfrac{20}{n}\left(\dfrac{3}{2}\right)^n
%\\
%&
+\dfrac{30}{n(n-1)}\left(\dfrac{9}{5}\right)^n-\dfrac{40}{n(n-1)}\left(\dfrac{3}{2}\right)^n+\dfrac{10}{n(n-1)} \\
&=\dfrac{5(19n+5)}{4n(n-1)}\left(\dfrac{9}{5}\right)^n-\dfrac{20(n+1)}{n(n-1)}\left(\dfrac{3}{2}\right)^n-\dfrac{5(3n-11)}{4n(n-1)}.
\end{align*}
In the last step, we have $|u_n| \leq (\frac{9}{5})^n$, as for $n>25$, the following two inequalities hold:
\begin{eqnarray}
\dfrac{5(19n+5)}{4n(n-1)} & \leq & 1 \nonumber \\
-\dfrac{20(n+1)}{n(n-1)}\left(\dfrac{3}{2}\right)^n-\dfrac{5(3n-11)}{4n(n-1)} & \leq & 0. \nonumber
\end{eqnarray}
Thus, the claim is proved. $\Box$

\medskip

We now consider the set $\mathcal{B}\equiv \{z\in\mathbb{C}\ :\ |z|\leq \frac{1}{2} \}$, and the partition $U(z)=\sum_{k=0}^\infty u_kz^k= U_1(z) + U_2(z)$, $U_1(z) \equiv \sum_{k=0}^{100}u_kz^k$ and $U_2(z) \equiv \sum_{k=101}^\infty u_kz^k$. Using the bound for $|u_n|$ from Lemma \ref{uuu}, for each $z \in \mathcal{B}$ we have
\begin{equation}
\label{upper-bound}
|U_2(z)|\leq \sum_{k=101}^\infty |u_k| \, |z|^k \leq \sum_{k=101}^\infty \left(\dfrac{9}{5}\right)^k \left(\dfrac{1}{2}\right)^k= 10\left(\dfrac{9}{10}\right)^{101} \approx 0.0002390525900.
\end{equation}
Next, we need a lower bound for $|U_1(z)|$.

\begin{lemm}
\label{lower-bound}
We have $\min_{z\in\partial\mathcal{B}}|U_1(z)|\geq \dfrac{3}{1000}$.
\end{lemm}
\noindent \emph{Proof.}
We obtain the result by considering a function
\[
G(t)\equiv \left[\sum_{k=0}^{100}u_k\cos(kt)\left(\dfrac{1}{2}\right)^k\right]^2+\left[\sum_{k=0}^{100}u_k\sin(kt)\left(\dfrac{1}{2}\right)^k\right]^2. \]
$G(t)$ has period $2 \pi$, with $G(\pi - t) = G(\pi + t)$, if $t \in [0,\pi]$. For $|z|\in \partial\mathcal{B}$ we can write $z = \frac{1}{2} [\cos t + i \sin t]$ for $t \in [0,2\pi)$, and thus
\begin{eqnarray*}
|U_1(z)| &=& \left|\sum_{k=0}^{100} u_k \left[\left( \frac{1}{2} \right) [\cos t + i \sin t] \right]^k \right| = \left|\sum_{k=0}^{100} u_k \cos(k t) \left( \frac{1}{2} \right)^k + i \sum_{k=0}^{100} u_k \sin(k t) \left( \frac{1}{2} \right)^k \right|
= \sqrt{G(t)}.
\end{eqnarray*}

By using the bound in Lemma \ref{uuu}, we have the following inequality
\begin{align}
|G'(t)|&=\left|2\left[\sum_{k=0}^{100}u_k\cos(kt)\left(\dfrac{1}{2}\right)^k\right]\left[-\sum_{k=0}^{100}ku_k\sin(kt)\left(\dfrac{1}{2}\right)^k\right]\right.\nonumber\\
 &\quad\qquad\left.+2\left[\sum_{k=0}^{100}u_k\sin(kt)\left(\dfrac{1}{2}\right)^k\right]\left[\sum_{k=0}^{100}ku_k\cos(kt)\left(\dfrac{1}{2}\right)^k\right]\right|\nonumber\\
&\leq 2\left| \sum_{k=0}^{100}u_k\cos(kt) \left( \dfrac{1}{2} \right)^k \right| \left| \sum_{k=0}^{100}ku_k\sin(kt)\left(\dfrac{1}{2}\right)^k \right| \nonumber\\
 &\quad\qquad +2\left| \sum_{k=0}^{100}u_k\sin(kt)\left(\dfrac{1}{2}\right)^k\right| \left|\sum_{k=0}^{100}ku_k\cos(kt)\left(\dfrac{1}{2}\right)^k\right|\nonumber\\
&\leq 2\left[ \sum_{k=0}^{100}|u_k| |\cos(kt)| \left( \dfrac{1}{2} \right)^k \right] \left[ \sum_{k=0}^{100}k|u_k| |\sin(kt)|\left(\dfrac{1}{2}\right)^k \right] \nonumber\\
 &\quad\qquad +2\left[ \sum_{k=0}^{100}|u_k| |\sin(kt)|\left(\dfrac{1}{2}\right)^k\right] \left[\sum_{k=0}^{100}k|u_k| |\cos(kt)|\left(\dfrac{1}{2}\right)^k\right]\nonumber\\
 &\leq 4\left[\sum_{k=0}^{100}\left(\dfrac{9}{10}\right)^k\right]\left[\sum_{k=0}^{100}k\left(\dfrac{9}{10}\right)^k\right] \approx 3598.862135.\label{bound-der}
\end{align}

We set $\mathcal{I}=\{\frac{k\pi}{1000000}\ :\ k\in\mathbb{Z},0\leq k\leq1000000 \}$. A numerical calculation shows that
\begin{equation}
\label{min-on-I}
\min_{t\in \mathcal{I}}G(t)=G(0) \approx 0.01949528529.
\end{equation}
With these preparations complete, we prove our claim by showing that
\begin{equation}
\label{suff-bound}
\displaystyle\min_{t\in[0,\pi]}G(t)\geq \dfrac{9}{1000000}.
\end{equation}
We prove Eq.~(\ref{suff-bound}) by contradiction. Suppose there exists $t_0\in[0,\pi]$ such that $G(t_0)<\frac{9}{1000000}$. Then we can find $t_1\in \mathcal{I}$ such that
\begin{equation}
\label{dis-t1-t0}
\vert t_1-t_0\vert\leq \dfrac{\pi}{2000000}.
\end{equation}
By the Mean Value Theorem, we can find $c \in (t_0,t_1)$ such that $G(t_1)-G(t_0)=G'(c)(t_1-t_0)$. From Eqs.~(\ref{bound-der}) and (\ref{dis-t1-t0}),
\begin{equation}
\label{eqUpper}
\dfrac{1800\pi}{1000000}\geq |G'(c)(t_1-t_0)|=|G(t_1)-G(t_0)|\geq G(t_1)-G(t_0).
\end{equation}
However, because $t_1\in{\mathcal I}$, by Eq.~(\ref{min-on-I}), we have
\[
G(t_1)-G(t_0)\geq G(0)-G(t_0)\geq\dfrac{1}{100}-\dfrac{9}{1000000}= \dfrac{9991}{1000000}.
\]
This result contradicts the upper bound in Eq.~(\ref{eqUpper}). Thus, Eq.~(\ref{suff-bound}) holds and the claim has been proven.
$\Box$

\medskip

Next, we study the root of $U_1(z)$ inside ${\mathcal B}$.

\begin{lemm}
\label{lem10}
The polynomial $U_1(z)$ has a unique (simple) root $\beta$ inside $\mathcal{B}$, with $\beta \approx 0.4889986317$.
\end{lemm}

\noindent \emph{Proof.}
First, by the Intermediate Value Theorem, there exists a real root $\beta$ with $0< \beta < \frac{1}{2}$, as we can numerically compute $U_1(0) \, U_1(\frac{1}{2})<0$ for the polynomial $U_1(z)$. Thus, we must prove
\[
\frac{U_1(z)}{z-\beta}=\frac{U_1(z)-U_1(\beta)}{z-\beta}=\sum_{k=0}^{100}u_k\frac{z^k-\beta^k}{z-\beta}=
\sum_{k=0}^{100}u_k\sum_{\ell=0}^{k-1}\beta^{k-1-\ell}z^{\ell}=\sum_{\ell=0}^{99}\left(\sum_{k=\ell+1}^{100}u_k\beta^{k-1-\ell}\right)z^{\ell}
\]
satisfies $\left| {U_1(z)}/(z-\beta) \right| > 0$ in $\mathcal{B}$.

To do so, we first use the bisection method for root-finding to numerically approximate $\beta$ by
\[
\tilde{\beta}=\dfrac{1101127027820569}{2251799813685248} \approx 0.4889986317,
\]
with the approximation error
\begin{equation}
\label{erralfa}
\vert\beta-\tilde{\beta}\vert\leq \frac{1}{2^{50}}.
\end{equation}

Then, we define the polynomial
$$Q(z)\equiv\sum_{\ell=0}^{99} a_{\ell}z^{\ell}, \text{ with } a_{\ell} \equiv \sum_{k=\ell+1}^{100}u_k\tilde{\beta}^{k-1-\ell},$$
through which we can write
\begin{eqnarray}
\frac{U_1(z)}{z-\beta} & = & Q(z)+(\beta-\tilde{\beta})R(z), \nonumber \\
R(z) & \equiv & \sum_{\ell=0}^{99}\left(\sum_{k=\ell+1}^{100}u_k\frac{\beta^{k-1-\ell}-\tilde{\beta}^{k-1-\ell}}{\beta-\tilde{\beta}}\right)z^{\ell}
=\sum_{\ell=0}^{99}\left(\sum_{k=\ell+2}^{100}u_k\sum_{j=0}^{k-2-\ell}\beta^j\tilde{\beta}^{k-2-\ell-j}\right)z^{\ell}. \nonumber
\end{eqnarray}
Note that on ${\mathcal B}$,
\begin{equation}
\label{upp-Rz}
\vert R(z)\vert\leq \sum_{\ell=0}^{99}\sum_{k=\ell+2}^{100}\sum_{j=0}^{k-2-\ell} |u_k| |\beta|^j |\tilde{\beta}|^{k-2-\ell-j} |z|^{\ell} \leq \sum_{\ell=0}^{99}\sum_{k=\ell+2}^{100}\sum_{j=0}^{k-2-\ell}\left(\frac{9}{5}\right)^k \left(\frac{1}{2}\right)^{k-2} \approx 3234.224489,
\end{equation}
where we used the bound for $|u_n|$ from Lemma \ref{uuu} and the fact that $\beta,\tilde{\beta},|z| \leq \frac{1}{2}$.

Next, let us consider the function
$$S(r,\theta)\equiv \sum_{\ell=0}^{99}a_{\ell}r^{\ell}\cos(\ell\theta)$$
defined over the rectangle $(r,\theta) \in [0, \frac{1}{2}] \times [0,\pi],$ where $S(r,\theta) = \Re(Q(z))$ if $z=r[\cos(\pm \theta) + i \sin(\pm \theta)] \in \mathcal{B}$.
We need the following bound for the gradient of $S$:
\begin{eqnarray}
\nonumber
| \nabla S | &=& \left| \left(\sum_{\ell=0}^{99} \ell a_{\ell} r^{\ell-1} \cos(\ell \theta),\sum_{\ell=0}^{99} -\ell a_{\ell} r^{\ell} \sin(\ell \theta) \right) \right| = \left| \sum_{\ell=0}^{99} \left( \ell a_{\ell} r^{\ell-1} \cos(\ell \theta), -\ell a_{\ell} r^{\ell} \sin(\ell \theta) \right) \right| \\\nonumber
&=& \left| \sum_{\ell=0}^{99} \ell a_{\ell} r^{\ell-1} \left( \cos(\ell \theta), - r \sin(\ell \theta) \right) \right| \leq \sum_{\ell=0}^{99} \ell |a_{\ell}| |r|^{\ell-1} |\left( \cos(\ell \theta), - r \sin(\ell \theta) \right)| \\\label{bound-grad}
&\leq & \sum_{\ell=0}^{99} \ell |a_{\ell}| |r|^{\ell-1} \leq \sum_{\ell=0}^{99} \ell |a_{\ell}| \left( \frac{1}{2} \right)^{\ell-1} \approx 89.628949.
\end{eqnarray}
Here, we have made use of $|r| < \frac{1}{2}$ and for $|r| < 1$, $\sqrt{\cos^2 x + r^2 \sin^2 x} \leq  \sqrt{\cos^2 x + \sin^2 x} = 1$.

%Next, for $z=r[\cos(\theta) + i \sin(\theta)] \in \mathcal{B}$, let $S(r,\theta):=\Re(Q(z))$, i.e.,
%\[
%S(r,\theta)=\sum_{\ell=0}^{99}a_{\ell}r^{\ell}\cos(\ell\theta).
%\]
%We need the following (rough) bound for the gradient of $S$:
%\begin{equation}\label{bound-grad}
%\vert\nabla S\vert \leq \sqrt{2}\sum_{\ell=1}^{99}\ell|a_{\ell}|\left(
%\dfrac{1}{2}\right)^{\ell-1}=126.75447\ldots.
%\end{equation}

A numerical calculation shows that over the grid $\mathcal{I} \equiv \{(\frac{k}{2000},\frac{j\pi}{1000})\ :\ (k,j)\in\mathbb{Z}^2,0\leq k,j\leq 1000 \}$, we have
\begin{equation}
\label{min-on-grid}
\min_{(r,\theta)\in\mathcal{I}}|S(r,\theta)|=\left|S\left(\dfrac{1}{2},\dfrac{502\pi}{1000}\right)\right| \approx 0.9518894218.
\end{equation}
We now show---with a similar method to that used to prove Lemma \ref{lower-bound}---that
\begin{equation}
\label{suff-bound-2}
\min_{(r,\theta)\in[0, \frac{1}{2}]\times[0,\pi]}\vert S(r,\theta)\vert\geq\frac{3235}{2^{50}}.
\end{equation}

Suppose for contradiction that there exists $z_0=(r_0,\theta_0) \in[0, \frac{1}{2}]\times[0,\pi]$ such that $| S(r_0,\theta_0) |<3235/ 2^{50}$. Then let us take $z_1=(r_1,\theta_1)\in\mathcal{I}$ such that
\begin{equation}
\label{dis-z1-z0}
|z_1-z_0|<\sqrt{\frac{1}{16}+\frac{\pi^2}{4}} \bigg(\frac{1}{1000}\bigg)\leq\frac{1}{500}.
\end{equation}
By the Mean Value Theorem, there exists a point $(r,\theta)$ on the line segment from $(r_0,\theta_0)$ to $(r_1,\theta_1)$ such that
\[
\nabla S(r,\theta)\cdot (z_1-z_0)=S(r_1,\theta_1)-S(r_0,\theta_0),
\]
where $\cdot$ is the inner product of $\mathbb{R}^2$. By using the Cauchy-Schwarz inequality together with (\ref{bound-grad}), (\ref{min-on-grid}) and (\ref{dis-z1-z0}), the assumption $| S(r_0,\theta_0) |<3235/2^{50}$ would thus give
\begin{eqnarray*}
%\dfrac{130}{500}
\frac{90}{500} &\geq & |\nabla S(r,\theta)||z_1-z_0|\geq|\nabla S(r,\theta)\cdot (z_1-z_0)|=|S(r_1,\theta_1)-S(r_0,\theta_0)| \\
&\geq & |S(r_1,\theta_1)|-|S(r_0,\theta_0)| \geq \dfrac{9}{10}-\dfrac{3235}{2^{50}} > 0.89,
\end{eqnarray*}
which is a contradiction. Hence, Eq.~(\ref{suff-bound-2}) holds.

Finally, because for $z\in{\mathcal B}$ we have
$$\vert Q(z) \vert \geq \vert\Re(Q(z))\vert \geq \min_{(r,\theta)\in[0, \frac{1}{2}]\times[0,\pi]}\vert S(r,\theta)\vert,$$
by using Eqs.~(\ref{erralfa}), (\ref{upp-Rz}), and (\ref{suff-bound-2}) it follows that in ${\mathcal B}$,
\begin{eqnarray*}
\left|\dfrac{U_1(z)}{z-\beta}\right|&=& \left|Q(z)+(\beta-\tilde{\beta})R(z)\right| \geq \left| |Q(z)| - |(\tilde{\beta}-\beta) R(z) | \right| \geq
\frac{3235}{2^{50}} - |(\tilde{\beta}-\beta)| | R(z) | \\
&\geq & \frac{3235}{2^{50}} - \dfrac{|R(z)|}{2^{50}} > \frac{3235}{2^{50}} - \dfrac{3234.224489\ldots}{2^{50}} > 0.
\end{eqnarray*}
This concludes the proof.
$\Box$

\medskip

Combining Lemmas \ref{lower-bound} and \ref{lem10} with the inequality in Eq.~(\ref{upper-bound}), we obtain the following proposition.
\begin{prop}
\label{rootOfU}
The function $U(z)$ has a unique (simple) root $\alpha$ inside $\mathcal{B}$, where $\alpha \approx 0.4889986317$.
\end{prop}

\noindent \emph{Proof.}
For the decomposition $U(z) = U_1(z) + U_2(z)$, Eq.~(\ref{upper-bound}) together with Lemma \ref{lower-bound} gives for $z \in \partial \mathcal{B}$
$$|U_1(z)| \geq \frac{3}{1000} > 0.00025 > |U_2(z)|.$$
Hence, from Rouch\'{e}'s theorem, inside ${\mathcal B}$ the function $U(z)$ has the same number of roots (considered with multiplicity) as polynomial $U_1(z)$. From Lemma \ref{lem10}, we know that $U_1(z)$ has one (simple) root inside ${\mathcal B}$.

The only remaining step is the numerical computation of $\alpha$, whose first ten digits turn out to coincide with the constant $\beta$ found in Lemma \ref{lem10} as the root of $U_1(z)$ inside $\mathcal{B}$. We again decompose $U(z)$:
$$U(z)=\sum_{k=0}^\infty u_kz^k=\sum_{k=0}^{500}u_kz^k+\sum_{k=501}^\infty u_kz^k = \tilde{U}_1(z)+\tilde{U}_2(z).$$
Note that from our bound for $|u_k|$ (Lemma \ref{uuu}), for each $z \in \mathcal{B}$ we have
\begin{equation}
\label{eqU2}
|\tilde{U}_2(z)| \leq \sum_{k=501}^{\infty} |u_k| \, |z|^k \leq \sum_{k=501}^\infty \left(\dfrac{9}{5}\right)^k \left(\dfrac{1}{2}\right)^k= 10\bigg(\frac{9}{10}\bigg)^{501} \leq 10^{-21}.
\end{equation}

Let us now consider
\begin{eqnarray}
\alpha ' & = & \dfrac{550563513910285}{1125899906842624}  \approx 0.48899863172938484723  \nonumber \\
\alpha'' & = & \dfrac{1101127027820571}{2251799813685248} \approx 0.48899863172938529132. \nonumber
\end{eqnarray}
These values were chosen using the bisection method such that
\[
\tilde{U}_1(\alpha')=2.708185805\ldots\cdot10^{-16}\quad\text{and}\quad\tilde{U}_1(\alpha'')=-4.953373282\ldots\cdot10^{-15}.
\]
From the bound of $|\tilde{U}_2(z)|$ in Eq.~(\ref{eqU2}), it is clear that $U(\alpha')>0$ and $U(\alpha'')<0.$ Let $\alpha$ be the unique root of $U(z)$ in $\mathcal{B}$, which by the Intermediate Value Theorem must be a real root in $(\alpha',\alpha'')$, and let $\epsilon \equiv \alpha-\alpha' \leq10^{-14}$. Note that
\[
\dfrac{1}{\alpha'}-\dfrac{1}{\alpha}=\dfrac{\epsilon}{\alpha'(\alpha'+\epsilon)}\leq\dfrac{\epsilon}{(\alpha')^2}\leq5\cdot10^{-14}.
\]
Thus, we can use
\begin{eqnarray}
\alpha'        & = & 0.48899863172938484723 \nonumber \\
(\alpha')^{-1} & = & 2.0449954971518340953 \nonumber
\end{eqnarray}
to approximate $\alpha$ and $\alpha^{-1}$, respectively.
$\Box$

%%%%%%%%%%%%%%%%%%%%%%%%%%%%%%%%%%%%%%%%%%%%%%%%%%%%%%%%%%%%%%
%%%%%%%%%%%%%%%%%%%%%%%%%%%%%%%%%%%%%%%%%%%%%%%%%%%%%%%%%%%%%%
%%%%%%%%%%%%%%%%%%%%%%%%%%%%%%%%%%%%%%%%%%%%%%%%%%%%%%%%%%%%%%
\medskip
\noindent
{\footnotesize
{\bf Acknowledgments} This work developed from discussions at the Banff International Research Station. Support was provided by a Rita Levi-Montalcini grant from the Ministero dell'Istruzione, dell'Universit\`a e della Ricerca (FD), grants MOST-104-2923-M-009-006-MY3 and MOST-107-2115-M-009-010-MY2 (MF, ARP), and National Institutes of Health grant R01 GM131404 (NAR). }
{\footnotesize
\bibliographystyle{chicago}
\bibliography{map3}}
%%%%%%%%%%%%%%%%%%%%%%%%%%%%%%%%%%%%%%%%%%%%%%%%%%%%%%%%%%%%%%
%%%%%%%%%%%%%%%%%%%%%%%%%%%%%%%%%%%%%%%%%%%%%%%%%%%%%%%%%%%%%%
%%%%%%%%%%%%%%%%%%%%%%%%%%%%%%%%%%%%%%%%%%%%%%%%%%%%%%%%%%%%%%

\end{document}